\documentclass[final,5p,times,twocolumn,authoryear]{elsarticle}

\usepackage{amssymb}
\usepackage{amsmath}
\usepackage{bm}
\usepackage{lipsum}
\usepackage{subcaption}
\usepackage{float}
\setcitestyle{square}

\usepackage{pgfplots}                   
\pgfplotsset{width=7cm, compat=newest}
\usetikzlibrary{patterns}
\tikzset{every label/.style={font=\footnotesize,inner sep=1pt}}
\newcommand{\stencilpt}[4][]{\node[circle,fill = red,draw,inner sep=1.5pt,label={above:#4},#1] at (#2) (#3) {}}

\usepackage{tikz}
\usepackage{hyperref}
\usepackage{ulem}
\usepackage{wrapfig}
\usepackage{standalone}
\usetikzlibrary{cd}
\usetikzlibrary{positioning}
\usetikzlibrary{calc}
\usetikzlibrary{fit}
\tikzset{set/.style={draw,circle,inner sep=0pt,align=center}}

\definecolor{morange}{RGB}{255,127,14}
\definecolor{mblue}{RGB}{31,119,180}
\definecolor{mred}{RGB}{214,39,40}
\definecolor{mpurple}{RGB}{148,103,189}
\definecolor{mgreen}{RGB}{44,160,44}

\definecolor{ccccccc}{RGB}{204,204,204}
\definecolor{cffffff}{RGB}{255,255,255}
\newcommand{\barnode}[4][]{\node[label={below:#4}] at (#2) (#3) {}}
\journal{Report article}

\begin{document}

\begin{frontmatter}

\title{A data-driven learned discretization approach in finite volume schemes for hyperbolic conservation laws and varying boundary conditions}

\author[first, second]{G. de Romèmont}
\author[first]{F .Renac}
\author[first]{J. Nunez}
\author[first]{D. Gueyffier}
\author[second]{F. Chinesta}
\affiliation[first]{organization={ONERA/DAAA},
            addressline={29 Av. de la Division Leclerc}, 
            city={Châtillon},
            postcode={92320}, 
            state={France}}

\affiliation[second]{organization={ENSAM},
            addressline={151 Bd de l'Hôpital}, 
            city={Paris},
            postcode={75013}, 
            state={France}}


\begin{abstract}
This paper presents a data-driven finite volume method for solving 1D and 2D hyperbolic partial differential equations . This work builds upon the prior research \citep{bar2019learning,zhuang2021learned,kochkov2021machine} incorporating a data-driven finite-difference approximation of smooth solutions of scalar conservation laws, where optimal coefficients of neural networks approximating space derivatives are learned based on accurate, but cumbersome solutions to these equations. We extend this approach to flux-limited finite volume schemes for hyperbolic scalar and systems of conservation laws. We also train the discretization to efficiently capture discontinuous solutions with shock and contact waves, as well as to the application of boundary conditions. The learning procedure of the data-driven model is extended through the definition of a new loss, paddings and adequate database. These new ingredients guarantee computational stability, preserve the accuracy of fine-grid solutions, and enhance overall performance. Numerical experiments using test cases from the literature in both one- and two-dimensional spaces demonstrate that the learned model accurately reproduces fine-grid results on very coarse meshes.
\end{abstract}


\end{frontmatter}




\section{Introduction}
\label{introduction}
Numerical fluid mechanics typically deals with nonlinear hyperbolic partial differential equations (PDEs) to represent complex fluid phenomena of various types. Regardless of the smoothness of initial or boundary data, discontinuities may develop after a finite time \cite[Sec. 2.4.2]{toro2000centred}. Therefore, appropriate numerical methods must be designed for these scenarios, particularly for approximating discontinuous solutions. The simulation of complex physical systems is fundamental, in practice, the sheer scale of these simulations—whether in terms of spatial resolution, temporal accuracy, or the complexity of the underlying physical processes—means that simulating systems at realistic levels is not feasible due to the computational time required for a given simulation.
\par 
The numerical integration of fluid equations on a coarse grid relative to a fully resolved integration of a given problem is a field of active research. Practically, coarse graining facilitates the computation of significantly larger systems. In this context, a number of works assessed the use of Neural Networks (NN) to solve PDEs using different methodologies have been developed. One enables neural ansatz to simultaneously align with target data while reducing a PDE residual \citep{lagaris1998artificial,raissi2018hidden, raissi2019physics, sun2020surrogate}. Other methodologies such as operator learning \citep{tran2021factorized, hao2023gnot, li2020fourier}, or purely learned surrogates \citep{pfaff2020learning, ronneberger2015u, dolz2018hyperdense, sanchez2020learning} have found application across a wide range of applied mathematics problems. 
\par
In practice, a significant advancement has been the capacity to precisely satisfy certain physical or numerical constraints by incorporating learned models into a fixed equation of motion. To this end and for forward problems, a number of hybrid physics-ML methods have emerged to train Neural Networks in the frame of existing numerical methods. This has been applied to a variety of topics such as turbulence \citep{list2022learned, stachenfeld2021learned, kochkov2021machine}, with learned artificial viscosity \citep{schwander2021controlling} or the use of constrained PINNs \citep{patel2022thermodynamically, coutinho2023physics, yu2018deep, kharazmi2019variational}. 
\par
Recent studies have employed machine learning (ML) techniques to enhance the resolution of shocks in hyperbolic PDEs, these problems are particularly challenging due to the potential for instabilities to occur. This has been achieved by utilizing the finite-volume method and several enhanced approaches have been proposed, among which popular ones are optimised shock-capturing schemes \citep{stevens2020enhancement, kossaczka2021enhanced}, learning flux-limiters \citep{nguyen2022machine,schwarz2023reinforcement}, learning corrections \citep{discacciati2020controlling} or coupling finite-volume with PINNs methodology \citep{ranade2021discretizationnet, cen2024deep}. Other ideas within the finite-volume framework have been explored, with LSTM model \citep{stevens2020finitenet} or for variational finite volume scheme \citep{zhou2024machine}. Other methods have also been developed on unstructured grids \citep{jessica2023finite, li2023finite}.

\par 
With the same objective in mind, \cite{bar2019learning} proposed a learned interpolation of the gradient that achieves the same accuracy as traditional finite difference methods but with much coarser grid resolution. The method has been extended into classical finite volume solvers to encompass passive scalar advection \citep{zhuang2021learned} and the Navier-Stokes equations \citep{kochkov2021machine}  . These methods are specific to the underlying equations and necessitate the training of a coarse-resolution solver with high-resolution ground truth simulations.
\par
The objective of this paper is to extend the approach introduced by Bar-Sinai to hyperbolic conservation laws, to greatly improve the treatment of shocks and discontinuities and adapt these method to various boundary conditions. The construction of a high-order ML-based numerical scheme is beset by a number of challenges, which must be addressed in order to ensure the scheme's efficacy :
\begin{itemize}
    \item \textit{Stability.} The presence of discontinuities presents a challenge in approximating the solution using high-order numerical methods, as the appearance of Gibbs oscillations \citep{wilbraham1848certain} can lead to instabilities, which could potentially leads the solution to grow without bound in a finite amount of time (diverge) or reduced robustness of any numerical scheme.
    \item \textit{Accuracy.}  The convergence rate of the overall numerical solution is conditioned by the accuracy of the discretization scheme in time and the evaluation of spatial derivatives. The accuracy of the spatial derivatives reconstruction is inextricably linked to the mesh size. This process is particularly challenging in the context of super-resolution.
    \item \textit{Computational performance.} The objective of the learned models is to produce high-fidelity simulations with a significant reduction in computational resources. 
\end{itemize}
\par 

The methodology proposed by \cite{bar2019learning} has been adapted to hyperbolic partial differential equations (PDEs), resulting in notable improvements in accuracy, robustness, and solution quality when solving conservation laws with neural networks. In section \ref{sec:model_problem}, we present an introduction to hyperbolic conservation laws with different PDE examples, the finite volume solver and boundary conditions are presented in section \ref{sec:finite_volume_solver}.
Along with the methodology, a number of regularisers have been introduced section \ref{sec:data_driven} to ensure that the neural network solution to the forward problems is physically achievable, particularly in situations involving strong shocks. Furthermore, the use of a limiter enables the training process to be completed without the need for unrolled steps, thereby significantly reducing the overall training time. 
The algorithm has been subjected to theoretical analysis section \ref{sec:toy_pb} with a view to demonstrating its potential for achieving gains in a linear equation. Furthermore, the method has been tested on a variety of equations and cases from the literature in section \ref{sec:results}, demonstrating excellent results in terms of accuracy and stability. We conduct numerical experiments to investigate computational performance of the algorithm in section \ref{sec:Computational}.
\par
The code for the Burgers 1D hyperbolic equation is available at \url{https://github.com/guigzair/Burgers\_1D}.

%
%
\section{Model problem}\label{sec:model_problem}

\subsection{Nonlinear hyperbolic conservation laws} 

We are here interested in the approximation of first-order nonlinear hyperbolic conservation laws and consider initial and boundary value problems in $d\geq1$ space dimensions of the form

\begin{subequations}\label{eq:hyperbolic}
\begin{align}
    \partial_t \textbf{w} + \nabla\cdot\textbf{f}(\textbf{w}) &= 0 \quad \text{in }\Omega\times(0,\infty), \label{eq:hyperbolic_a}\\
		\textbf{w}(\bm{x},0) &= \textbf{w}_0(\bm{x}) \quad \text{in }\Omega,  \label{eq:hyperbolic_b}\\
    \bm{B}(\textbf{w},\textbf{w}_{bc},\bm{n}) &= 0 \quad \text{on }\partial\Omega\times(0,\infty),  \label{eq:hyperbolic_c}
\end{align}
\label{eq:hyperbolic_eq}
\end{subequations}

\noindent with $\Omega \subset \mathbb{R}^d$ a bounded domain, $\textbf{w}:\Omega \times [0,T] \rightarrow \Omega^a\subset\mathbb{R}^r$ denotes the vector of $r$ conservative variables with initial data $\textbf{w}_0\in L^\infty(\Omega)$ and $\bm{f}:\Omega^a  \rightarrow \mathbb{R}^r \times \mathbb{R}^d$ are the physical fluxes. The solution is known to lie within a convex set of admissible states $\Omega^a$. Boundary conditions are imposed on $\partial\Omega$ through the boundary operator in (\ref{eq:hyperbolic_c}) and some prescribed boundary data $\textbf{w}_{bc}$ defined on $\partial\Omega$, while $\bm{n}$ denotes the unit outward normal vector to $\partial\Omega$. The operator $\bm{B}$ depends on the type on condition to be imposed and the equations under consideration. Examples are provided in section \ref{sec:BC_examples}.

Solutions to (\ref{eq:hyperbolic}) may develop discontinuities in finite time even if $\textbf{w}_0$ is smooth, therefore the equations have to be understood in the sense of distributions. Weak solutions are not necessarily unique and (\ref{eq:hyperbolic}) must be supplemented with further admissibility conditions to select the physical solution. We here focus on entropy inequalities of the form 
\cite[sec. I.5]{godlewski2013numerical}\cite[sec. III]{dafermos2005hyperbolic}\cite[sec. III.8]{leveque1992numerical}

\begin{equation}\label{eq:entropy_ineq}
    \partial_t \eta(\textbf{w}) + \nabla\cdot\bm{q}(\textbf{w})\leq 0,
\end{equation}

\noindent where $\eta:\Omega^a\rightarrow\mathbb{R}$ is a convex entropy function and $\bm{q}(\textbf{w}):\Omega^a\rightarrow\mathbb{R}^d$ is the entropy flux satisfying the compatibility condition 

\begin{equation}
    \nabla_{\textbf{w}} \eta(\textbf{w})^T \times \nabla_{\textbf{w}}\textbf{f}_i(\textbf{w}) = \nabla_{\textbf{w}}\bm{q}_i(\textbf{w})^T, \quad 1\leq i\leq d.
\end{equation}

The inequality (\ref{eq:entropy_ineq}) becomes an equality for smooth solutions, while strict convexity of $\eta(\textbf{w})$ implies hyperbolicity of (\ref{eq:hyperbolic_a}) \citep{godlewski2013numerical,dafermos2005hyperbolic}.

We here aim at satisfying these properties at the discrete level with our data-driven finite volume scheme.

\subsection{Model examples} \label{model_examples}
\subsubsection{Scalar equations}

We will first investigate the linear stability properties of the data-driven scheme with a von Neumann analysis. To this end, we consider the linear scalar advection equation in one space dimension:

\begin{equation}\label{eq:lin_advection}
    \partial_t w + a\partial_xw = 0,
\end{equation}
with $a>0$ the transport velocity which is assumed to be constant. 
\par

We will consider numerical experiments with the nonlinear scalar Burgers' equation in one space dimension:

\begin{equation}\label{eq:Burgers_1D}
    \partial_t w + \partial_x\Big(\frac{w^2}{2}\Big) = 0,
\end{equation}

\noindent for which the physical weak solutions satisfy the entropy inequality (\ref{eq:entropy_ineq}) for the entropy pair $\eta(w)=\tfrac{w^2}{2}$ and $q(w)=\tfrac{w^3}{3}$, together with a maximum principle on $w$ which may define $\Omega^a=[\min_\Omega w_0,\max_\Omega w_0]$. Solutions to scalar equations are known to be of bounded total variation \citep{godlewski2013numerical,dafermos2005hyperbolic}.

\subsubsection{Compressible Euler equations}

We will also consider numerical experiments with the compressible Euler equations in one, $d=1$, and two space dimensions, $d=2$, for which

\begin{align*}
    \textbf{w} &= (\rho, \rho \bm{v}^T, E)^T, \\
    \bm{f}(\textbf{w}) &= \big(\rho \bm{v},\rho \bm{v}\bm{v}^T + p{\bf I}_d, (E+p)\bm{v}\big)^T,
\end{align*}

\noindent where ${\bf I}_d$ is the identity matrix of size $d$, while $\rho$, $\bm{v}$ and $E$ denote the density, velocity vector and total energy, respectively. We close the system with an equation of state for a polytropic ideal gas law, so  $E=\tfrac{p}{\gamma - 1} + \tfrac{1}{2}\rho |\bm{v}|^2$ with $\gamma>1$ the ratio of specific heats and $p$ the pressure. Note that it will be convenient to consider the numerical discretization as a function of the primitive variables 

\begin{equation*}
    \textbf{u} = (\rho, \bm{v}^T, p)^T.
\end{equation*}
Note that for scalar equations, $u = w$.
\par
This system is hyperbolic in every unit direction $\bm{n}$ over the set of admissible states $\Omega^a=\{\textbf{w}\in\mathbb{R}^{d+2}:\,\rho>0,\bm{v}\in\mathbb{R}^d,E-\tfrac{1}{2}\rho |\bm{v}|^2>0\}$ with eigenvalues $\bm{v}\cdot\bm{n}-c$, $\bm{v}\cdot\bm{n}$, $\bm{v}\cdot\bm{n}+ c$, where $c=\big(\tfrac{\gamma p}{\rho}\big)^{1/2}$ is the speed of sound. An entropy pair is given by

\begin{equation}
    \eta = -\rho s, \quad
    \bm{q}(\textbf{w}) = -\rho s \bm{v},
\end{equation}

\noindent where $s = C_v \log\left(\tfrac{P}{\rho^\gamma}\right)$ is the physical specific entropy defined by the Gibbs relation.

\subsection{Boundary conditions}\label{sec:BC_examples}

To a given boundary $\Gamma\subset\partial\Omega$ in (\ref{eq:hyperbolic_c}), we associate an admissible boundary state

\begin{equation}\label{eq:boundary_state}
 \textbf{w}_\Gamma:\Omega^a\times\mathbb{S}^{d-1}\ni(\textbf{w},\bm{n})\mapsto\textbf{w}_\Gamma(\textbf{w},\bm{n})\in\Omega^a,
\end{equation}

\noindent which we assume to be admissible, $\textbf{w}_\Gamma(\Omega^a\times\mathbb{S}^{d-1})\subset\Omega^a$, and consistent: ${\bf B}(\textbf{w},\textbf{w}_{bc},\bm{n})=0$ implies $\textbf{w}_\Gamma(\textbf{w},\bm{n})=\textbf{w}$.


%
\paragraph{Boundary operators}

Several boundary conditions have been implemented. The treatment of boundary conditions is a critical aspect in every CFD solver, as the way boundaries are handled can significantly influence the overall accuracy and stability of the solution, which is essential for capturing the true behavior of the system

In the above scenario, location $a$ is outside the domain and corresponds to a ghost cell, $\Gamma$ is the boundary interface, and location $d$ is inside the physical domain. The unit normal vector $\bm{n} = [n_x, n_y]^T$ points outwards of the domain. 
\par 
The flux Jacobian $\bm{A}(\textbf{w},\bm{n}) = \partial_{\textbf{w}} \bm{f}(\textbf{w})\cdot\bm{n}$ has eigenvalues $(\lambda_i)_{1\leq i\leq r}$ and is diagonalized as $A(\textbf{w},  \bm{n}) = R^{-1} \Lambda R$ with $\Lambda = \text{diag}(\lambda_i)$. Depending on the sign of the eigenvalues, we introduce the upwind decomposition of $A(\textbf{w},  \bm{n})$ as follows
\begin{equation*}
    A(\textbf{w},  \bm{n}) = A^+(\textbf{w},  \bm{n}) + A^-(\textbf{w},  \bm{n}), 
\end{equation*}

\noindent with 
\begin{equation*}    
    \quad A^\pm(\textbf{w},  \bm{n}) = R^{-1} \Lambda^\pm  R, \quad \Lambda^\pm=\dfrac{1}{2}\text{diag}\big(\lambda_i\pm|\lambda_i|\big).
\end{equation*}
\par 
Each eigenvalue is associated to a given characteristics. The characteristics are used to determine the number of boundary conditions to be imposed \citep{goncalvesdasilva:cel-00556980,hartmann2002adaptive}, which correspond to the inflow characteristics that propagate from outside to inside the computational domain.
\par 

\paragraph{Compressible Euler equations}

The flux Jacobian $\bm{A}(\textbf{w},\bm{n}) = \partial_{\textbf{w}} \bm{f}(\textbf{w})\cdot\bm{n}$ has eigenvalues

\begin{equation*}
    \lambda_1=\bm{v}\cdot\bm{n}-c,\quad \lambda_2=\dots=\lambda_{d+1}=\bm{v}.\bm{n}, \quad \lambda_{d+2}=\bm{v}\cdot\bm{n} + c.
\end{equation*}

\paragraph{Farfield conditions} We can apply characteristic boundary conditions on the farfield boundary from a freestream state $\textbf{w}_{\infty}$ with 
\begin{equation}
    B(\textbf{w}, \textbf{w}_{bc}, \bm{n}) = A^- (\textbf{w} - \textbf{w}_{\infty}),
\end{equation}
\noindent where $\textbf{w}_{\infty}$ denotes a given freestream state. For instance, for a supersonic inflow boundary condition, we have $A^-=A$, corresponding to the Dirichlet boundary conditions
\begin{equation}
    \textbf{w}_{\Gamma}(\textbf{w}, \bm{n}) =  \textbf{w}_{\infty},
\end{equation}
while for the supersonic outflow boundary condition, we have $A^+=A$ and $A^-=0$, corresponding to an extrapolation condition
\begin{equation}
    \textbf{w}_{\Gamma}(\textbf{w}, \bm{n}) =  \textbf{w}.
\end{equation}

\paragraph{Slip boundary condition} Considering an impermeability condition, $\bm{v}.\bm{n} = 0$, at a wall $\Gamma_w \subset \partial \Gamma$, the associated boundary data is $\textbf{w}_\Gamma(\textbf{w}, \bm{n}) = (\rho, \rho (\bm{v}- (\bm{v}.\bm{n})\bm{n})^T, \rho E)^T$. The condition is commonly imposed  through the use of a mirror state $2 \textbf{w}_\Gamma(\textbf{w}, \bm{n}) - \textbf{w} = (\rho, \rho (\bm{v}- 2(\bm{v}.\bm{n})\bm{n})^T, \rho E)^T$. The mirror state $\bm{u}^+_\Gamma(\bm{u}^-, \bm{n})$ follows from imposing a linear reconstruction interpolating the right and left states at the interface : $\frac{1}{2}(\textbf{w}^- + \textbf{w}^+_\Gamma(\textbf{w}^-, \bm{n})) = \textbf{w}_\Gamma(\textbf{w}^-, \bm{n})$ hence $\textbf{w}^+_\Gamma(\textbf{w}^-, \bm{n})) = 2 \textbf{w}^-_\Gamma(\textbf{w}, \bm{n}) - \textbf{w}^-$

\paragraph{Periodic condition} 
This boundary is established when physical geometry of interest and expected flow pattern are of a periodically repeating nature. This reduces computational effort in our problems.

%
%
\section{Finite volume solver} \label{sec:finite_volume_solver}

We describe below the finite volume solver used to generate the data. We use a MUSCL reconstruction of the slopes \citep{van1979towards} on 1D and 2D Cartesian meshes to get a formally second-order scheme. In section \ref{finite_volume_solver}, we introduce the scheme in 1D for the sake of clarity, while the treatment of boundary conditions is described in section \ref{sec:FV_BCs}.

\subsection{Finite volume solver and reference solution} \label{finite_volume_solver}

The domain is partitioned with a uniform and Cartesian mesh with cells $I_i = [x_{i-\frac{1}{2}}, x_{i+\frac{1}{2}}]$ of size $\Delta x_i$ in 1D and we look for an approximate solution to problem (\ref{eq:hyperbolic}) of the form of a piecewise constant solution where the degrees of freedom approximate the cell-averaged solution in each cell:

\begin{equation*}
    \overline{\textbf{w}}_{i}(t) \simeq \dfrac{1}{\Delta x_i}\int_{I_i} \textbf{w}(\bm{x},t) dV. 
\end{equation*}

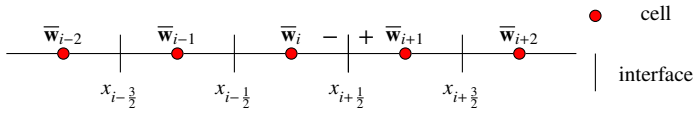
\begin{figure}[H]
    \begin{center}
        \begin{tikzpicture}
          \stencilpt{-3,0}{i-2}{$\overline{\textbf{w}}_{i-2}$};
          \stencilpt{-1.5,0}{i-1}{$\overline{\textbf{w}}_{i-1}$};
          \stencilpt{ 0,0}{i}  {$\overline{\textbf{w}}_{i}$};
          \stencilpt{ 1.5,0}{i+1}{$\overline{\textbf{w}}_{i+1}$};
          \stencilpt{ 3,0}{i+2}{$\overline{\textbf{w}}_{i+2}$};
          \draw (i-2) -- (i-1)
                (i-1) -- (i)
                (i)   -- (i+1)
                (i+1) -- (i+2);
            \draw (-0.75,-0.25) -- (-0.75,0.25);
            \barnode{ -0.75,-0.25}{}{$x_{i-\frac{1}{2}}$};
            \draw (-2.25,-0.25) -- (-2.25,0.25);
            \barnode{ -2.25,-0.25}{}{$x_{i-\frac{3}{2}}$};
            \draw (0.75,-0.25) -- (0.75,0.25);
            \draw (0.75,-0.25) -- (0.75,0.25);
            \draw (0.75,0.00) node[above left]  {$-$};
            \draw (0.75,0.00) node[above right] {$+$};
                \barnode{ 0.75,-0.25}{}{$x_{i+\frac{1}{2}}$};
            \draw (2.25,-0.25) -- (2.25,0.25);
                \barnode{2.25,-0.25}{}{$x_{i+\frac{3}{2}}$};
            \draw (i-2)-- (-3.75,0);
            \draw (i+2) -- (3.75,0);
            \stencilpt{ 4.,0.5}{label}{};
            \barnode{4.8,0.8}{label}{cell};
            \draw (4,-0.5) -- (4.,0.);
            \barnode{4.8,0}{label}{interface};
        \end{tikzpicture}
    \caption{Example of a 1D finite volume mesh.}
    \label{fig:enter-label} 
    \end{center}
\end{figure}

Integrating (\ref{eq:hyperbolic_a}) over a fixed volume $V$ and using the divergence theorem over the surface $\partial V$ gives
\begin{equation}
    \partial_t \int_{V} \textbf{w}(t)dV + \oint_{\partial V}\textbf{f}(\textbf{w})\cdot \textbf{n}dS = 0
\end{equation}
The discrete scheme in cell $I_i$ with an explicit Euler time integration approximates this conservation law as follows:
\begin{equation}
    \overline{\textbf{w}}_{i}^{n+1}-\overline{\textbf{w}}_{i}^{n}+\dfrac{\Delta t_n}{\Delta x_i}  \big(\hat{\textbf{f}}(\textbf{w}^{-,n}_{i+1/2}, \textbf{w}^{+,n}_{i+1/2}) -\hat{\textbf{f}}(\textbf{w}^{-,n}_{i-1/2}, \textbf{w}^{+,n}_{i-1/2})\big) = 0
\end{equation}

\noindent where $\overline{\textbf{w}}_{i}^{n}=\overline{\textbf{w}}_{i}(t^{n})$ and $\textbf{w}^{\pm,n}_{i+1/2}=\textbf{w}^{\pm}_{i+1/2}(t^n)$ with $t^{n+1}-t^n=\Delta t_n>0$ the time step, while $\hat{\textbf{f}}(\cdot, \cdot)$ denotes a two-point numerical flux that is assumed to be  consistent $\hat{\textbf{f}}(\textbf{w},\textbf{w})=\textbf{f}(\textbf{w})$. By $\textbf{w}^\mp_{i+1/2}$ we denote approximate values of the reconstructed left and right traces at interface $x_{i+1/2}$ (see Figure \ref{fig:enter-label}). 

As a reference solver, we use a second-order MUSCL \citep{van1979towards} finite volume scheme with a Rusanov numerical flux \citep{lax2005weak, rusanov1961calculation} and a van Albada limiter \citep{van1982comparative} which is used because it is a smooth differentiable function as opposed for example to the minmod slope limiter \citep{roe1986characteristic}. The flux can be written
\begin{equation*}
    \hat{\textbf{f}}(\textbf{w}^-, \textbf{w}^+) = \frac{1}{2}
    (\textbf{f}(\textbf{w}^-)+\textbf{f}(\textbf{w}^+))-
    \frac{1}{2}
    \rho(\textbf{w}^-,\textbf{w}^+)
    (\textbf{w}^+-\textbf{w}^-)
\end{equation*}

\noindent where $\rho(\textbf{w}^-,\textbf{w}^+)=\max|\nabla_{\bf w}{\bf f}(\textbf{w}^\pm)|$ denotes the maximum spectral radius of the Jacobian of the flux $\textbf{f}(\cdot)$. Defining $\Phi(r)=\tfrac{r^2+r}{r^2+1}$ as the van Albada limiter, the second-order MUSCL reconstruction is defined on the primitive variables $\textbf{u}$ as 
\begin{equation}
    \begin{cases}
        u^-_{i+1/2} = u_i + \frac{1}{2}\phi_i(u_{i+1}-u_i),\\
        u^+_{i+1/2} = u_{i+1} - \frac{1}{2}\phi_{i+1}(u_{i+2}-u_{i+1}), \\
        \phi_i = \Phi\left(\dfrac{u_{i}-u_{i-1}}{u_{i+1}-u_i}\right).
    \end{cases}
    \label{eq:MUSCL_reconstruction}
\end{equation}

\noindent for every component $u$ in $\textbf{u}$. The reference solution is computed at a sufficiently high resolution, then down-sampled keeping only every $R$-th sample to produce the test and training datasets.
\par 
The rational for using a first order time integration is as follows. Tests have been performed with a second order time discretization scheme using the Heun Runge-Kutta method but no improvement was observed in the produced reference solutions. It also increased the training time. 

\subsection{Boundary conditions}\label{sec:FV_BCs}
 The implementation of far-field boundary conditions for a flow problem is contingent upon two prerequisites. Primarily, the truncation of the domain must not have a discernible impact on the flow solution when compared to that of an infinite domain. Secondly, any outgoing disturbances must not be reflected back into the flow field. Frequently, the far-field values are not known; only the freestream values are typically available.
 \par
 In order to achieve this, the computational domain is modified by the introduction of so-called "ghost cells" situated outside of the boundary. The introduction of fictitious flow in the ghost cells will yield the desired boundary conditions at the edge.

\begin{figure}[H]
    \centering
        \begin{tikzpicture}
        \draw  (2,0.5) -- (2,3);
        \draw[dashed] (1,2) -- (2,2);
        \draw[] (2,2) -- (3,2) -- (4,2);
         \draw[dashed] (1,1.75) -- (1,2.25);
        \stencilpt{ 1.5,2}{i+2}{$\textbf{w}_{-1}$};
        \filldraw[black] (2,3) node[anchor=north east]{$\Gamma$};
        \draw (3,1.75) -- (3,2.25);
        \stencilpt{ 2.5,2}{i+2}{$\textbf{w}_{1}$};
        \draw (4,1.75) -- (4,2.25);
        \stencilpt{ 3.5,2}{i+2}{$\textbf{w}_{2}$};
        \draw[->] (2,1) -- (1,1) node[anchor=north west]{$\bm{n}$};
    \end{tikzpicture}
    \caption{Ghost cell approach at boundaries, the ghost cell is represented in dashed line}
    \label{fig:boundary}
\end{figure}
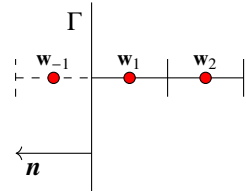

\paragraph{Supersonic inflow}
The conservative variables on the boundary are determined exclusively on the basis of freestream values. 
\begin{equation*}
    \textbf{w}_{\Gamma} = \textbf{w}_{-1}
\end{equation*}
\paragraph{Supersonic outflow} The conservative variables on the boundary are determined exclusively on the basis of inner flow field values. 
\begin{equation*}
    \textbf{w}_{\Gamma} = \textbf{w}_{1}
\end{equation*}
\paragraph{Subsonic inflow} 
Four characteristic variables are prescribed based on the freestream values, which are defined as follows: One characteristic variable is derived from the interior of the physical domain. This results in the following set of boundary conditions:
\begin{align*}
    p_{\Gamma} &= \frac{1}{2}\left[ p_{-1}+p_{1}-\rho_{0}c_{0}\bigl[(\bm{v}_{-1}- \bm{v}_{1})\cdot \bm{n}\bigr] \right] \\
    \rho_{\Gamma} &= \rho_{-1}+\frac{(p_{\Gamma}-p_{-1})}{c_{0}^{2}} \\
    \bm{v}_{\Gamma} &= \bm{v}_{-1}-\frac{\bm{n}(p_{-1}-p_{\Gamma})}{p_{0}c_{0}} 
\end{align*}
where $p_0$ and $c_0$ represent a reference state. The reference state is normally set equal to the state at the interior.
\paragraph{Subsonic outflow}
\begin{align*}
        p_{\Gamma} &= p_{-1}\\
        \rho_{\Gamma}&=\rho_{1}+\frac{\left(p_{\Gamma}-p_{1}\right)}{c_{0}^{2}}\\
        \bm{v}_{\Gamma}&=\bm{v}_{1}-\frac{\bm{n}\left(p_{1}-p_{\Gamma}\right)}{\rho_{0}c_{0}}
\end{align*}
with $p_{-1}$ being the prescribed static pressure. 
\paragraph{No-slip boundary} The implementation of the no-slip boundary condition can be simplified through the incorporation of ghost cells.
\begin{align*}
        p_{-1} &= p_{1} \\
        \rho_{-1} &= \rho_{1} \\
        \bm{v}_{-1} &= \bm{v}_{1} -  2(\bm{v}_{1}\cdot \bm{n})\bm{n}
\end{align*}
\paragraph{Periodic condition} 
Due to the periodicity condition, the first ghost-cell layer $\textbf{w}_{-1}$ corresponds to the  inner flow field value at the opposite periodic boundary.

%
%
\section{Data driven solution for hyperbolic equations} \label{sec:data_driven}

\subsection{Learning the derivatives}
One of our main goal is to improve accuracy without compromising generalization on a standard second-order finite volume solver as described in section \ref{finite_volume_solver}. To that extent, we consider a ML modeling to approach the derivatives. We expect that ML models can improve the accuracy of a CFD solver when  run  on  inexpensive-to-simulate  coarse  grids even in cases where the solutions exhibit shocks, discontinuities, or large gradients. By using a limiter and a MUSCL scheme, our solution is naturally limited to remove spurious oscillations \citep{harten1997high} that would otherwise occur around discontinuities. This ensures the robustness of the solution.
\par 
We assume a Cartesian grid in $d$ space dimensions. The accuracy of finite volume schemes relies on local high-order reconstructions of the solution, which are based on the evaluation of the space derivatives of the solution from its values in neighboring elements using finite differences. For a field $\textbf{w} (\textbf{x},t)$, the $i$th space derivateive can be approximated by
\begin{equation}
    \dfrac{\partial \textbf{w}}{\partial x_i}(\textbf{x},t) \approx \sum_{\substack{j=-k}}^k \alpha_j \textbf{w}(x_1, x_2, \dots, x_i + jh, \dots, x_N,t),
    \label{eq:gradient_conv}
\end{equation}
where ${x_1, ..., x_N}$ are the spatial coordinates and $h$ the size of the element which is here assumed constant. As an example, a first order approximation of the derivative on a field $w$ in 1D gives $\tfrac{\partial w}{\partial x} = \tfrac{w_{i+1}-w_{i}}{\Delta x}+\mathcal{O}(\Delta x)$ with $\alpha_1 = -\alpha_0 = \tfrac{1}{\Delta x}$ and $\alpha_{-1}=0$.
\par 
In this work, we follow the idea of \cite{bar2019learning} and \cite{zhuang2021learned} and extend it to finite volume schemes for hyperbolic conservation laws. The algorithm works as follows, at each time step, the neural networks generates a latent vector $\bm{\alpha} = (\alpha_{-k}, ...,\alpha_k)$ at each cell location for each variable based on the current primitive variables ${\bf u}$ of the PDE. The latent vectors are then used to reconstruct the space gradient of the primitives variables $\nabla \textbf{u}$ using a convolution as in equation (\ref{eq:gradient_conv}). All derivatives in equation (\ref{eq:MUSCL_reconstruction}) are then replaced by their respective reconstructions $\hat{\partial}_{x_j} \textbf{u}_{i}$ to predict the next time step. 
\par The new MUSCL scheme reads
\begin{equation}
    \begin{cases}
        u^-_{i+1/2} = u_i + \frac{1}{2}\phi_i\hat{\partial}u_i,\\
        u^+_{i+1/2} = u_{i+1} - \frac{1}{2}\phi_{i+1}\hat{\partial}u_{i+1}, \\
        \phi_i = \Phi\left(\dfrac{\hat{\partial}u_{i-1}}{\hat{\partial}u_{i}}\right).
    \end{cases}
    \label{eq:MUSCL_reconstruction_ML}
\end{equation}
\par
Because the derivative operator $\partial^n_{x_i} \cdot$ applied on a function $f$ at a point $x\in\Omega$ is local, independent of the size of the mesh and invariant by translation; the neural network used for learning the operator $\hat{\partial}\cdot$ also needs to verify these properties. 
The convolution operator satisfies the same properties as the derivative thus convolutionnal Neural Networks (CNNs) are used for the model. 
\par 
The training of a neural network inside a classic numerical solver is made possible by writing the entire program in a differentiable programming framework. Indeed the computation of the gradients of the loss is made possible using a differentiable CFD solver. This differentiable framework allows Automatic Differentiation \citep{baydin2018automatic} for any parameter of the solver. Multiple frameworks have been recently developed like TensorFlow \citep{abadi2016tensorflow}, Pytorch \citep{paszke2017automatic}, Flux.jl \citep{innes2018flux}, JAX \citep{bradbury2018jax} and are more and more efficient on hardware accelerator (GPUs, TPUs). These user-friendly frameworks make it easier to incorporate neural networks techniques into scientific codes. We have implemented our solver using TensorFlow Eager \citep{agrawal2019tensorflow}.

\subsection{Learned interpolation model}
The learned interpolation model is introduced in \cite{bar2019learning} and extended in \citep{zhuang2021learned, kochkov2021machine,  alieva2023toward}. It is a data-driven method which uses the outputs of the neural network to generate interpolation based on local variables or features of the PDE. 
\par
\begin{figure}[h]
    \vspace{-1em}
    \centering
    \begin{tikzpicture}[module/.style={draw, very thick, rounded corners, minimum width=15ex},
			normmodule/.style={module, fill=morange!40},
			convnmodule/.style={module, fill=mblue!40},
			paddingmodule/.style={module, fill=mpurple!40},
			arrow/.style={-stealth, thick, rounded corners},
		]
		\node (input) {Primitives $\textbf{u}$};
		\node[above of=input, normmodule, align=center, yshift=.0cm] (norm) {Normalisation};
		\node[above of=norm, paddingmodule, align=center] (pad) {Padding};
		\node[above of=pad, convnmodule, align=center, yshift=.1cm] (conv) {Convolutionnal \\ layer};
		\node[above of=conv, normmodule, align=center, yshift=.2cm] (cons) {Linear constraint};
		\node[above of=cons] (output) {$\bm{\alpha}$};
	
		\draw[arrow] (input) -- (norm);
		\draw[arrow] (norm) -- (pad); 
		\draw[arrow] (pad) -- (conv);
		\draw[arrow] (conv) -- (cons); 
		\draw[arrow] (cons) -- (output);
	
		\coordinate (attresidual) at ($(norm.south)!0.5!(input.north)$);
		\coordinate (ffnnresidual) at ($(conv.south)!0.5!(cons.north)$);
	
		\node[fit=(pad)(conv),draw, ultra thick, rounded corners, label=left:$\mathrm{N\times}$] (encoder) {};
\end{tikzpicture}
\caption{Architecture of our model. First, all primitives are normalized between 0 and 1. Subsequently, there are N blocks, comprising a padding layer and a convolutional layer. The padding layer serves to address the boundaries and to guarantee that the size of the output solution vector is identical to that of the input solution vector, as a convolution operation results in a reduction in vector size. Finally, the linear constraint guarantees the consistency of the model and reads $\sum_{j} \alpha_j=0$.}
\label{fig:architecture}
\vspace{1em}
\end{figure}
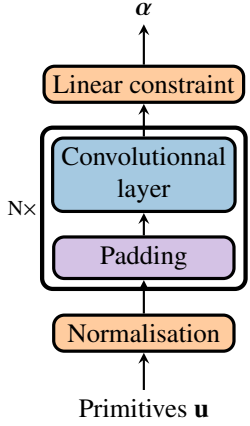

 
As shown on Figure \ref{fig:architecture}, the model is composed of several convolutionnal and padding layers  stacked together so that the output shape matches the input shape. The discrete convolution operator reduces the dimension of any input vector. Consequently, the output solution vector size is equal to the input solution vector size, as required by the padding layers. Paddings are adapted to the type of boundary condition during training and inference stages.
The normalization layer outputs all variables to $[-1, 1]$ with a min-max normalization. The linear constraint reads $\sum_{j=-k}^k \alpha_j=0$. It is similar to the ones introduced by \cite{bar2019learning} and \cite{zhuang2021learned} and enforces the derivative to be at least first-order accurate. (see next section below).

\subsection{Padding}
An essential feature of our model lies in how we handle the boundaries. It is imperative to guarantee that the model remains invariant by translation in space. Consequently, the convolutional layers illustrated in Figure \ref{fig:architecture} cannot be truncated in the direction normal to the boundary. In contrast, the flow variables are extended beyond the boundary via the use of a padding to set ghost cells. The extension method must ensure the preservation of the boundary condition on the coarse grid. The boundary condition is also applied as a padding for the periodic, inflow, or outflow boundary condition; however, it must be adapted for the slip boundary condition. \\

For the slip boundary condition, a first guess would be to mirror the flow across the wall, although this would result in the creation of an artificial discontinuity (see section \ref{sec:FV_BCs}) at the boundary. In order to guarantee the stability of the algorithm, the ML contribution of the learned discretization is set to zero at the interface for slip boundary conditions. This guarantees that the boundary condition is clearly defined and that the ML scheme is stabilized globally.

\subsection{Loss}
The training time tunes the ML parameters so that the weights and biases reduce the difference between a costly high-resolution simulation and a simulation generated by the model on a coarse grid. This refinement is achieved through supervised training where we use as the loss function $\mathcal{L}_{tot}$ the cumulative pointwise error between predicted and reference primitive variables. Because the numerical scheme is an inherently stable flux-limited second order scheme, it is unnecessary to accumulate the error over multiple steps for further stabilization.
\begin{equation}\label{eq:total_loss}
    \mathcal{L}_{tot} =  \left(|| \textbf{u}^{ref} - \hat{\textbf{u}}||_1+\lambda_2|| \textbf{u}^{ref} - \hat{\textbf{u}}||_2\right) + \mathcal{L}_{penalization},
\end{equation}
\par
where $\hat{\textbf{u}}$ is the machine learning solution, while $\textbf{u}^{ref}$ denotes the reference solution obtained on fine grids using the finite volume scheme. We have found that using a mix of the $L_1$ and $L_2$ (with a mean reduction operator) norms added with a penalization described later stabilizes efficiently the solution over time. Here the reference solution is a solution computed on a fine grid and projected on a coarser grid, the neural network solution is computed on the coarser grid. During training, the ratio $R = \frac{\Delta x_{coarse}}{\Delta x_{fine}}$ is 2 for stability reasons, but during inference, the ratio can be changed and increased.
\par 
Even if the method implemented without penalizations is stable and quite robust, penalization allows for an extra regularization and robustness, especially around shocks where spurious oscillations can appear even in the presence of a flux limiter. Extra penalizations terms are called soft constraints \citep{raissi2018hidden}. The Entropy and Total Variation Diminishing (TVD) penalizations were first introduced in \cite{patel2022thermodynamically} for hyperbolic PDEs. Here, we use
\begin{equation}
    \mathcal{L}_{penalization} = \lambda_{ent} \mathcal{L}_{ent} + \lambda_{TVD} \mathcal{L}_{TVD} + \lambda_{reg} \mathcal{L}_{reg}
\end{equation}

For the total loss (\ref{eq:total_loss}), the parameter $\lambda_2$ is set equal to 1, the three parameters $\lambda_{ent}$, $\lambda_{TVD} $ and $\lambda_{reg} $ are tuned through Bayesian optimization \citep{frazier2018tutorial} on the validation loss which is $\mathcal{L}_{tot} - \mathcal{L}_{penalization}$. 

\paragraph{Entropy inequality penalization}

Defining as in \cite{patel2022thermodynamically}, an entropy-flux loss can be written following equation (\ref{eq:entropy_ineq}) and be rewritten as a conservation law 
\begin{equation}
    \int_{\Omega \times T} \partial_t \eta + \nabla.\bm{q}d\omega dt \leq 0
\end{equation}
thus defining the inequality as a penalization,
\begin{equation}
    \mathcal{L}_{ent} = \sum_{\substack{j=1}}^{N_c}max\left(0,\int_{I_j \times T_j} \partial_t \eta + \nabla.\bm{q}d\omega dt \right)^2
\end{equation}
Where $I_j \times T_j \in \Omega \times T$ are the $N_c$ space-time cells on which we want to compute the loss. Giving in discrete form 
\begin{align*}
    K_j &= \int_{I_j \times T_j} \partial_t \eta + \nabla.\bm{q}d\omega dt \\
    K_j &= h^d(\eta(\bm{x_j}, t + dt) - \eta(\bm{x_j}, t)) + dt h^{d-1}\displaystyle \sum_{i_d}^d q_{i_d}(\bm{x_j},t) - q_{i_d}(\bm{x_j} - h\bm{e_{j}},t)
\end{align*}
where $\bm{q}\cdot\bm{e_{i_d}} = q_{i_d}$.
\paragraph{Total Variation penalization}
The TVD property on a 1D Cartesian mesh can often be given as 
\begin{equation}
    \forall t\in[O,T], TV(\textbf{u}^{t+\Delta t}) \leq TV(\textbf{u}^t),
\end{equation}
with
\begin{equation}
    T V\left(\textbf{u}^t\right)=\displaystyle \sum_{i_d}^d \sum_j^{N_{i_d}} |\textbf{u}( \dots, x_{i_d} + jh, \dots,t) - \textbf{u}( \dots, x_{i_d} + (j-1)h, \dots,t)|
\end{equation}
with $N_{i_d}$ the number of cells in the $i_d-th$ dimension. 
\par
Thus the penalization term can be expressed as
\begin{equation}
    \mathcal{L}_{TVD} = max(0,TV(\textbf{u}^{t+\Delta t}) - TV(\textbf{u}^t)).
\end{equation}
This TVD property can be applied to scalar conservation laws, with the existence of a total-variation bound. This is a highly desirable feature in numerical solutions of scalar conservation laws \citep{godlewski2013numerical}. Even if the existence of such bound does not hold in general in multiple space dimensions, nor for hyperbolic systems \citep{dafermos2005hyperbolic,rauch1986bv}, a soft TV penalization remains a popular approach to damp spurious oscillations \cite[Sec. 13.6]{toro2000centred}\cite{patel2022thermodynamically}. Furthermore, such penalization is applied as a soft constraint in the loss and we do not deal with the inequality directly.

\paragraph{Regularization term}
The regularization term is a penalization on the weights and biases of the neural network in order to avoid overfitting by preventing high values of the weights and biases. This has an effect of reducing oscillations for strong shocks. The penalization is written as:
\begin{equation}
    \mathcal{L}_{reg} = \lVert W \rVert_1
\end{equation}
with $W$ being all the parameters of the neural network.

\subsection{Dataset}
Multiple datasets with different boundary conditions have been used. The datasets are obtained on a fine grid and then projected on a coarser grid for training and for achieving super-resolution properties finite volume scheme. In particular, coarse-grained data are derived from finer data by employing the nodal values intrinsic to the finer resolution.
\par 
The training data is composed of sine wave and randomized piecewise constant initial conditions using pre-defined parametrized functions (for further information, see sec. \ref{sec:dataset}). In 2D, a circle in the center of a piecewise constant initial condition has been added in order to better learn derivatives in every direction. Special attention has been given to the quality of the database and the quality of the model as the results depend heavily on the  database used for the learning process. Here specifically, sine wave data is used for further regularization of the solution and the piecewise constant initial conditions are used to better learn the approximation of discontuitites in the solution.

\subsection{Model training} \label{sec:model_training}
\par
As opposed to the training of the algorithm in \cite{bar2019learning}, there is no need for rolling steps during training time to increase the scheme stability. Indeed, using a gradient reconstruction with a limiter allows for the solution to be more stable. 
\par
A transfer-learning method is used. The model is first trained on a periodic boundary dataset, then fine tuned on a slip boundary condition dataset. During the second phase of training, the weights $\alpha$ are set to zero on the edges in order to maintain the translation invariance property of the model. Periodic boundary conditions solutions have a tendency to contain stronger and more shocks. It makes the method more robust and fine tuning on slip boundary condition makes the method more generalizable due to the constraint imposed on the coefficients $\alpha$ on the edges. Without this training scheme, the ML method has a tendency to diverge with no-slip boundary conditions. 
\par
In practice, we use 3 blocks of convolutionnal layers and padding. Each convolutionnal layer has 32 filters with a kernel size of 3 (3x3 in 2D). The last layer has 4 filters (16 in 2D) in order to match the linear constraint and to output 3 coefficients for the reconstruction at each location on the mesh. The Adam optimizer \citep{da2014method} and SELU activation function are used. In 1D, the database is solved on 512 cells then projected on 256 cells for reference data, in 2D the database is solved on 256x256 cells then projected on a 128x128 grid for reference data. The training takes approximately 6h on a NVIDIA V100 GPU for the 2D case. 
%
%
\section{von Neumann stability analysis}\label{sec:toy_pb}
We here analyze the linear stability properties of the data driven scheme by performing a von Neumann stability analysis. We consider the scalar linear advection equation (\ref{eq:lin_advection}) with periodic boundary conditions. For a given real wavenumber $k$, (\ref{eq:lin_advection}) admits exact solution of the form $e^{i(kx-\omega t)}$, with $\omega=ka$ the complex frequency, as elementary solutions.

As described in section \ref{sec:finite_volume_solver}, the discrete scheme reads on a Cartesian mesh

\begin{equation*}
    w^{n+1}_j = w^{n}_j - C_o (w^{n}_{j} - w^{n}_{j-1}) - \dfrac{C_o \Delta x }{2}(\phi_j \hat{\partial}w^{n}_{j} - \phi_{j-1} \hat{\partial}w^{n}_{j-1})
\end{equation*}
with $C_o = \dfrac{a \Delta x}{\Delta t}$ the CFL number. In order to facilitate the analysis, we set the limiter to 1 and linearize the non-linear operator $\hat{\partial} \cdot$. This operator can be written as 
\begin{equation*}
    \hat{\partial}w_{j} = \dfrac{\alpha^{j}_{-1}w_{j-1} + \alpha^{j}_{0}w_{j}+ \alpha^{j}_{1}w_{j+1}}{\Delta x}
\end{equation*}
with $\alpha^{j}_{-1}, \alpha^{j}_{0} , \alpha^{j}_{1}$ being three non-linear functions depending on a localized stencil. We consider these three functions depend only on $(w_{j-1},w_{j},w_{j+1})$ for each position $j$. Considering $w_{j} = \Bar{w}+ w'_{j}$ ($\Bar{w}$ being the base state and $w'_{j}$ its perturbation at cell $j$), it gives at first order: 
\begin{equation*}
        \hat{\partial}w_{j} \simeq \dfrac{\Bar{\alpha}_{-1}w'_{j-1} + \Bar{\alpha}_{0}w'_{j}+ \Bar{\alpha}_{1}w'_{j+1}}{\Delta x} = \dfrac{\Bar{\alpha}_{-1}w_{j-1} + \Bar{\alpha}_{0}w_{j}+ \Bar{\alpha}_{1}w_{j+1}}{\Delta x}
\end{equation*}
With $\Bar{\alpha}_{k} = \alpha_{k}(\Bar{w}, \Bar{w}, \Bar{w})$ and where we have used the fact that $\Bar{\alpha}_{-1}+\Bar{\alpha}_{0}+\Bar{\alpha}_{1}=0$. Assuming periodic boundary conditions, we consider Fourier modes of the form $\hat{w}^{n}_j(k) = e^{i(kj\Delta x - \omega n\Delta t )}$ with $k$ being the given real wave number. We thus obtain 
\begin{align}
        \hat{w}^{n+1}_j(k) =& \hat{S}_{\Delta x}(k,\Delta t) \hat{w}^{n}_j(k) \\
        \hat{S}_{\Delta x}(k,\Delta t) =& 1-C_o (1-e^{-ik\Delta x})\Big(1 \nonumber\\ &+\dfrac{1}{2}\big(\bar{\alpha}_{-1}(e^{-ik\Delta x} - 1) + \bar{\alpha}_{1}(e^{ik\Delta x} - 1)\big)\Big)
        \label{eq:relation}
\end{align}

\begin{figure}[H]
     \centering
     \begin{subfigure}[b]{0.35\textwidth}
         \centering
         \includegraphics[width=\textwidth]{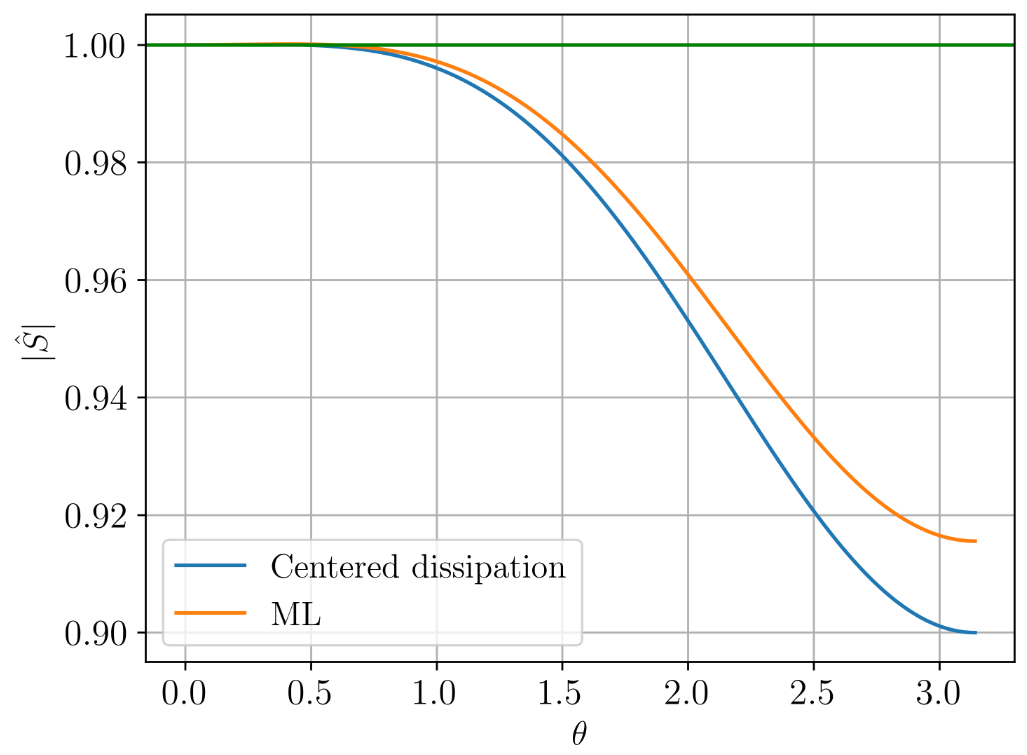}
         \caption{ $|\hat{S}_{\Delta x}(k,dt)|$}
         \label{fig:norm_centered_10}
     \end{subfigure}
     \hfill
     \begin{subfigure}[b]{0.35\textwidth}
         \centering
         \includegraphics[width=\textwidth]{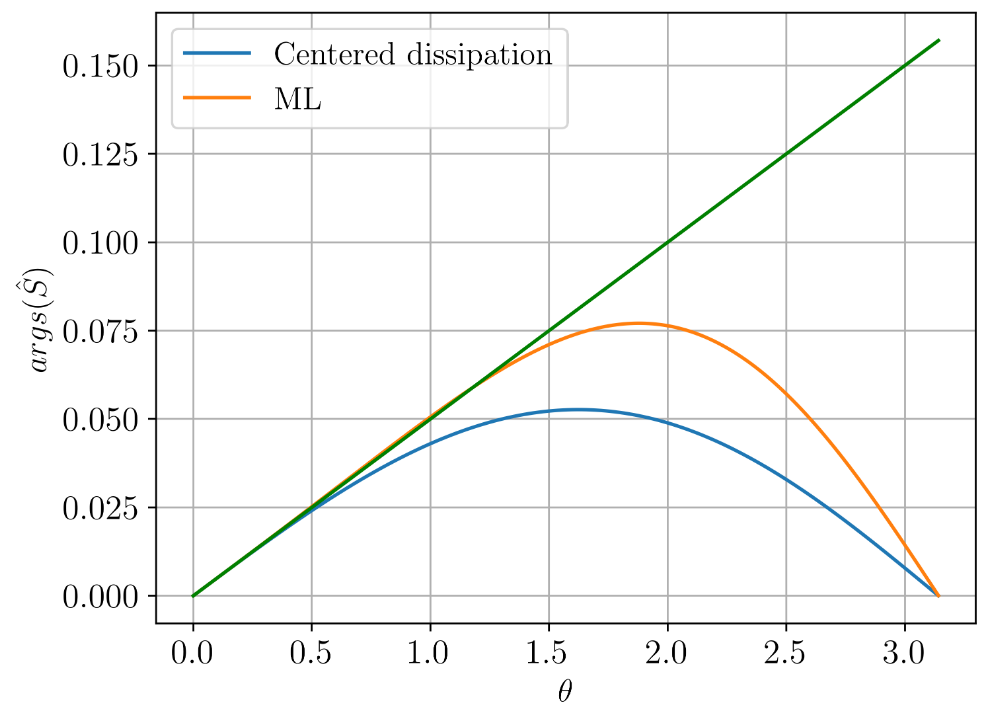}
         \caption{ $args(\hat{S}_{\Delta x}(k,dt))$}
         \label{fig:phase_centered_10}
     \end{subfigure}
     \caption{Numerical dissipation (a) and dispersion (b) of the ML scheme and the centered scheme. The green line corresponds to the exact dispersion relation for (\ref{eq:lin_advection}).}
\label{fig:scheme_analysis}
\end{figure}
Figure \ref{fig:scheme_analysis} shows the numerical dissipation and dispersion for eq. \ref{eq:relation} showing the ML scheme and a centered dissipation (taking $(\alpha_{-1}, \alpha_{0}, \alpha_{1}) = \left(-\frac{1}{2}, 0, \frac{1}{2}\right)$).

The figure shows that at first order the consistency and dispersion error is minimized for larger wavenumbers. 
%
%

\section{Results} \label{sec:results}
For all cases, an L2 norm is used for the error. All cases presented in this section are not present in the database.

\subsection{1D Burgers equation}
\begin{figure}[H]
     \centering
     \begin{subfigure}[b]{0.3\textwidth}
         \centering
         \includegraphics[width=\textwidth]{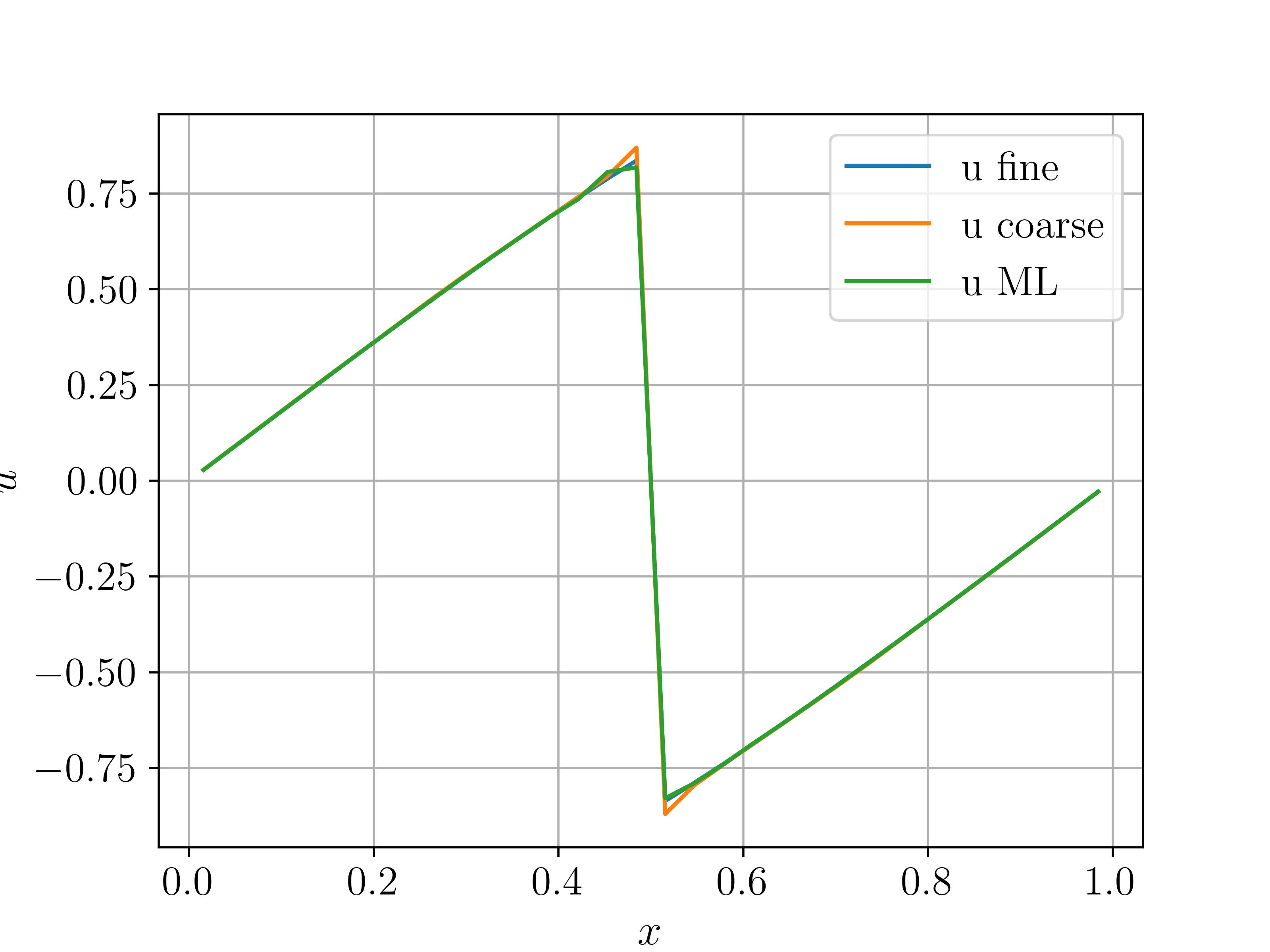}
         \caption{}
         \label{fig:sod_rho}
     \end{subfigure}
     \hfill
     \begin{subfigure}[b]{0.3\textwidth}
         \centering
         \includegraphics[width=\textwidth]{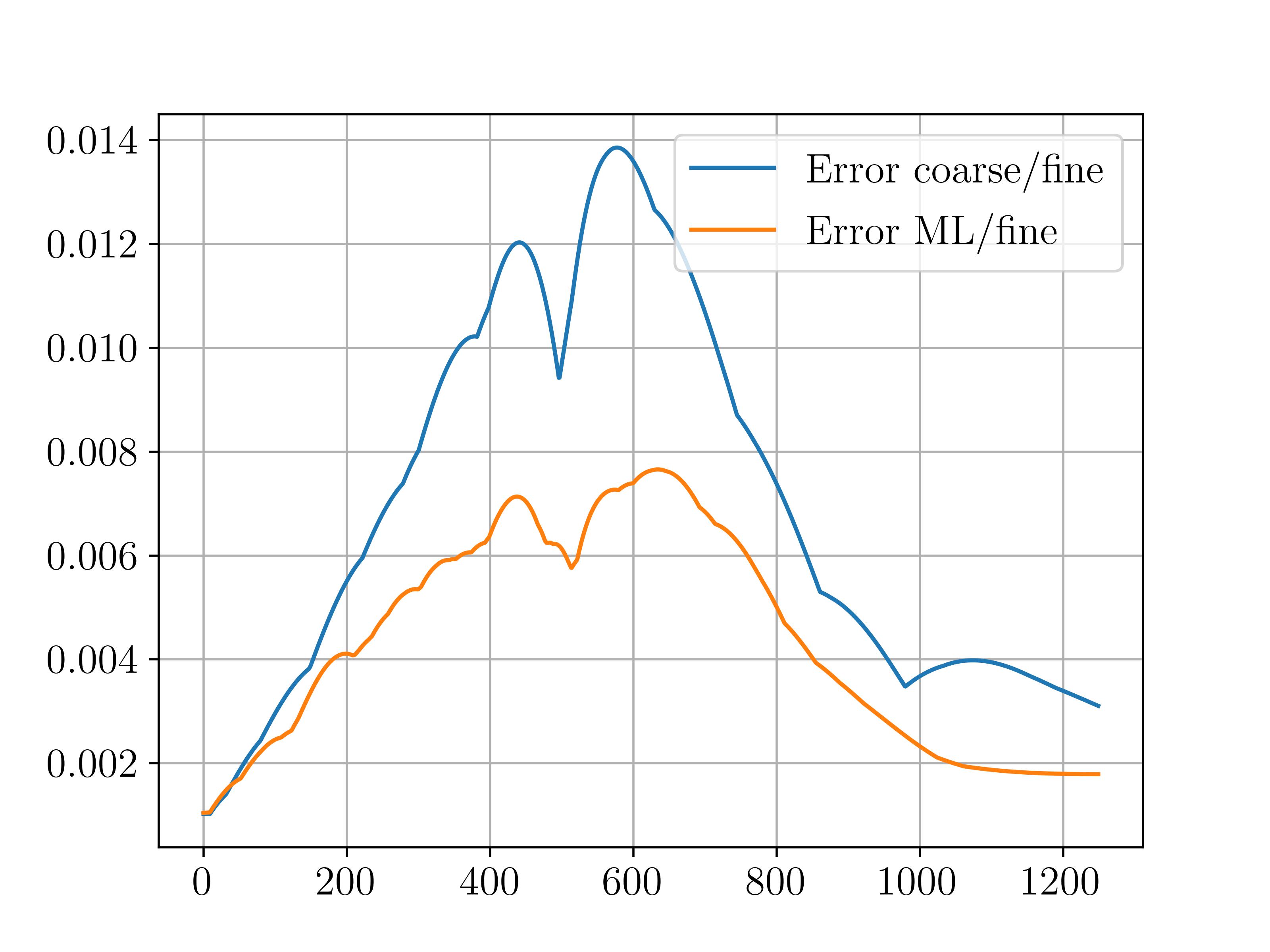}
         \caption{}
         
     \end{subfigure}
     \caption{Sine-wave test case, 1024 cells for fine integration, 32 cells for ML and
coarse solution, T = 0.39}
\label{fig:burgers_sin}
\end{figure}

\begin{figure}[H]
     \centering
     \centering
     \includegraphics[width=0.35\textwidth]{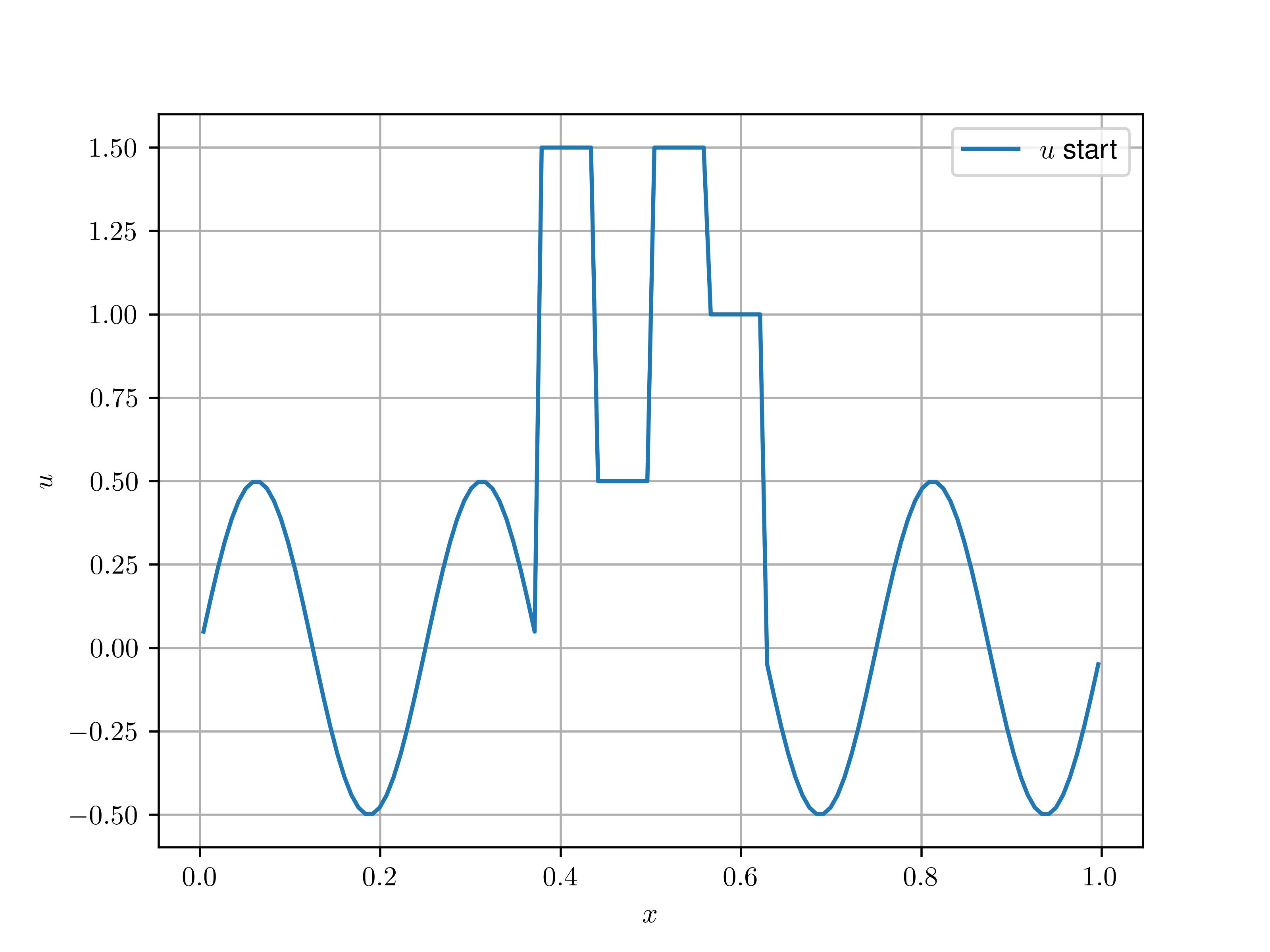}
     \label{fig:complex_initial}
     \caption{Complex test case initial condition}
\label{fig:burgers_complex_initial}
\end{figure}

The sine-wave test-case is the numerical solution of the 1D Burgers equation from a sinusoidal initial condition. The results for the reference solver and the ML one are represented Figure \ref{fig:burgers_sin} and show good results for this reference test case.
The initial condition \citep{ray2018artificial} is depicted in Figure \ref{fig:burgers_complex_initial} but is adapted to fit a $[0,1]$ domain and represent a more complex test case. 
The initial condition is defined as :
\begin{equation*}
    u_0(x)=
    \begin{cases}
        
        3 & \text { if } 3/8 \leq x \leq 3.5/8 \text { or } 4/8<x \leq 4.5/8 \\ 
        1 & \text { if }3.5/8<x<4/8, \\
        2 & \text { if } 4.5/8<x \leq 5/8, \\
        \sin (8\pi x) & \text { elsewhere. }  
    \end{cases}
\end{equation*}

\begin{figure}[H]
     \centering
     \begin{subfigure}[b]{0.3\textwidth}
         \centering
         \includegraphics[width=\textwidth]{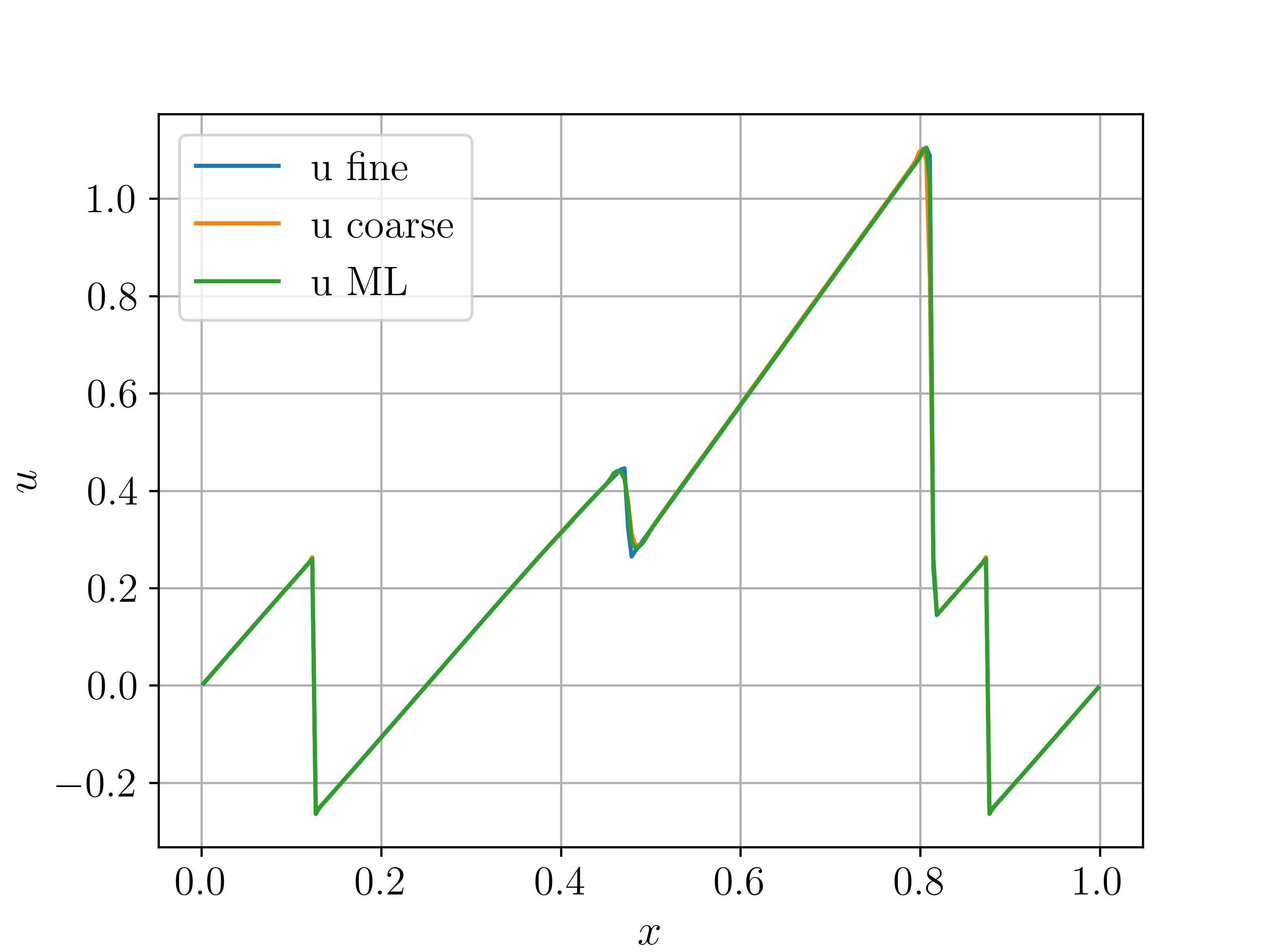}
         \caption{}
         \label{fig:sod_rho}
     \end{subfigure}
     \hfill
     \begin{subfigure}[b]{0.3\textwidth}
         \centering
         \includegraphics[width=\textwidth]{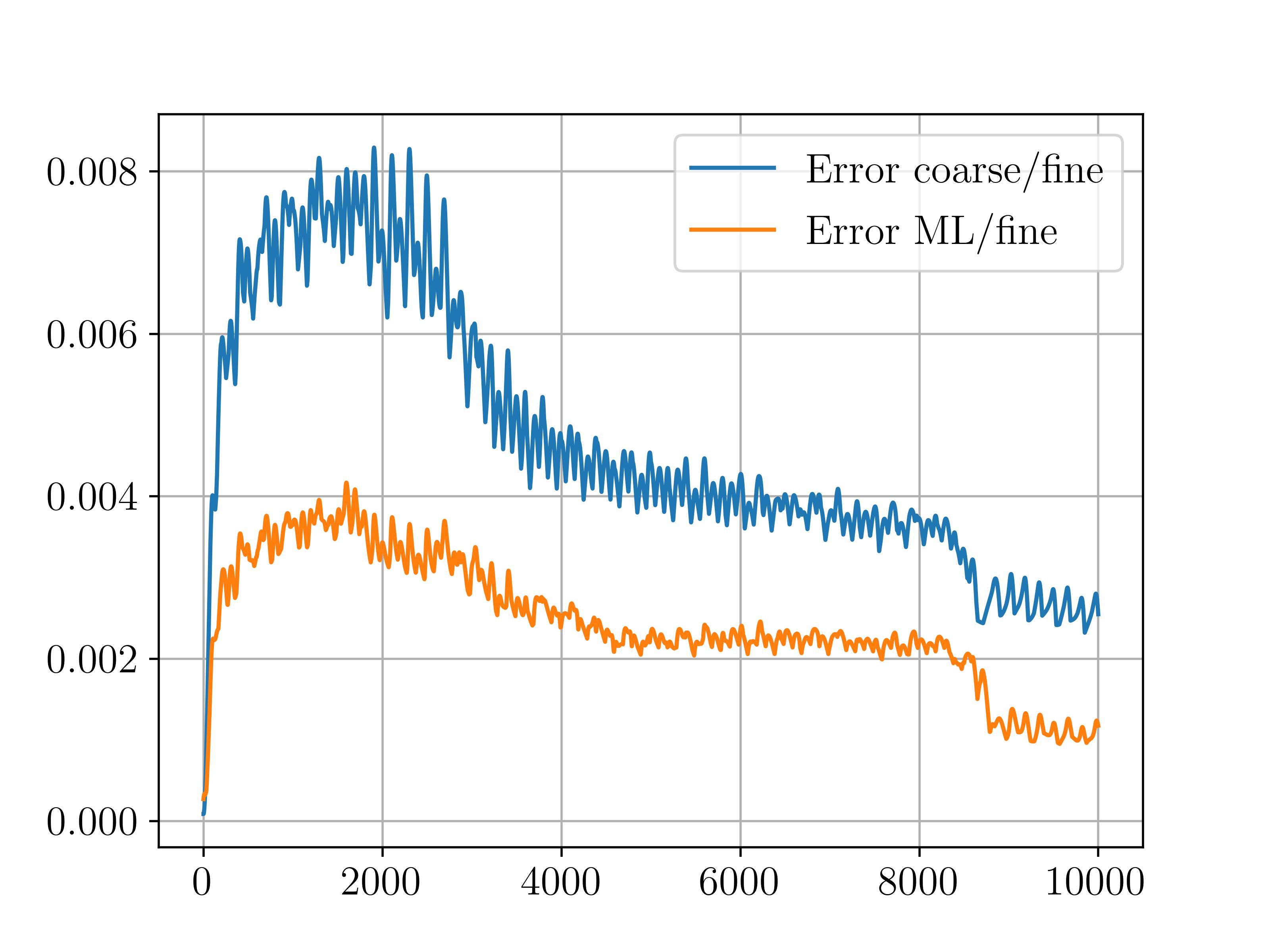}
         \caption{}
         
     \end{subfigure}
     \caption{Complex test case, 1024 cells for fine integration, 256 cells for ML and
coarse solution, T = 0.39}
\label{fig:burgers_complex}
\end{figure}
As seen on Figures \ref{fig:burgers_sin} and \ref{fig:burgers_complex}, the ML solver has a better precision mainly around discontinuities with sharp change of gradient where the limiter and its non-linearity has a strong effect.

\subsection{1D Euler equations}

\begin{figure}[H]
     \centering
     \begin{subfigure}[b]{0.3\textwidth}
         \centering
         \includegraphics[width=\textwidth]{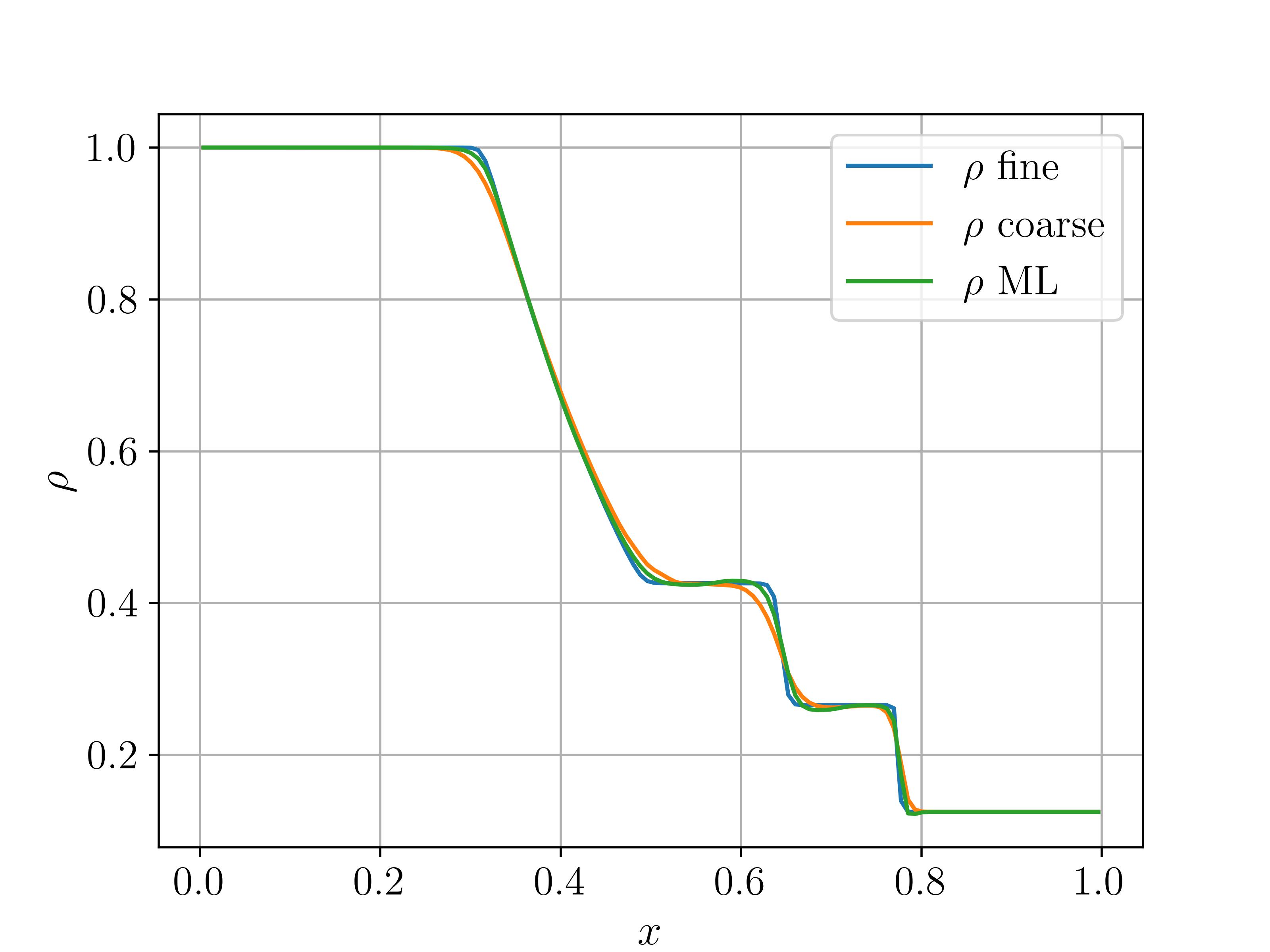}
         \caption{}
         \label{fig:sod_rho}
     \end{subfigure}
     \hfill
     \begin{subfigure}[b]{0.3\textwidth}
         \centering
         \includegraphics[width=\textwidth]{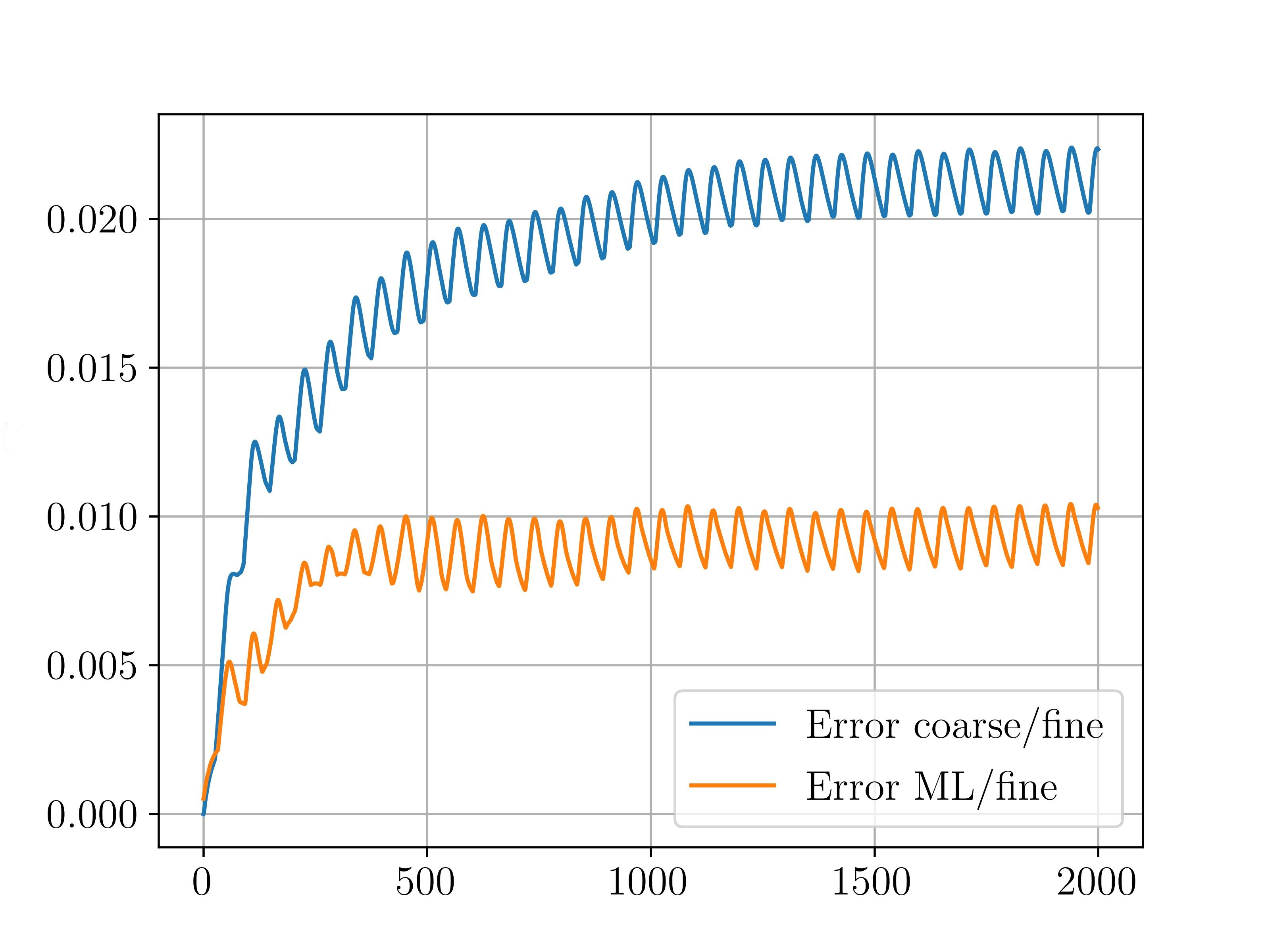}
         \caption{}
         
     \end{subfigure}
     \caption{Sod test case. 1024 cells for fine integration, 128 cells for ML and for coarse solutions, T =0.156.}
    \label{fig:sod}
\end{figure}

\begin{figure}[H]
     \centering
     \begin{subfigure}[b]{0.3\textwidth}
         \centering
         \includegraphics[width=\textwidth]{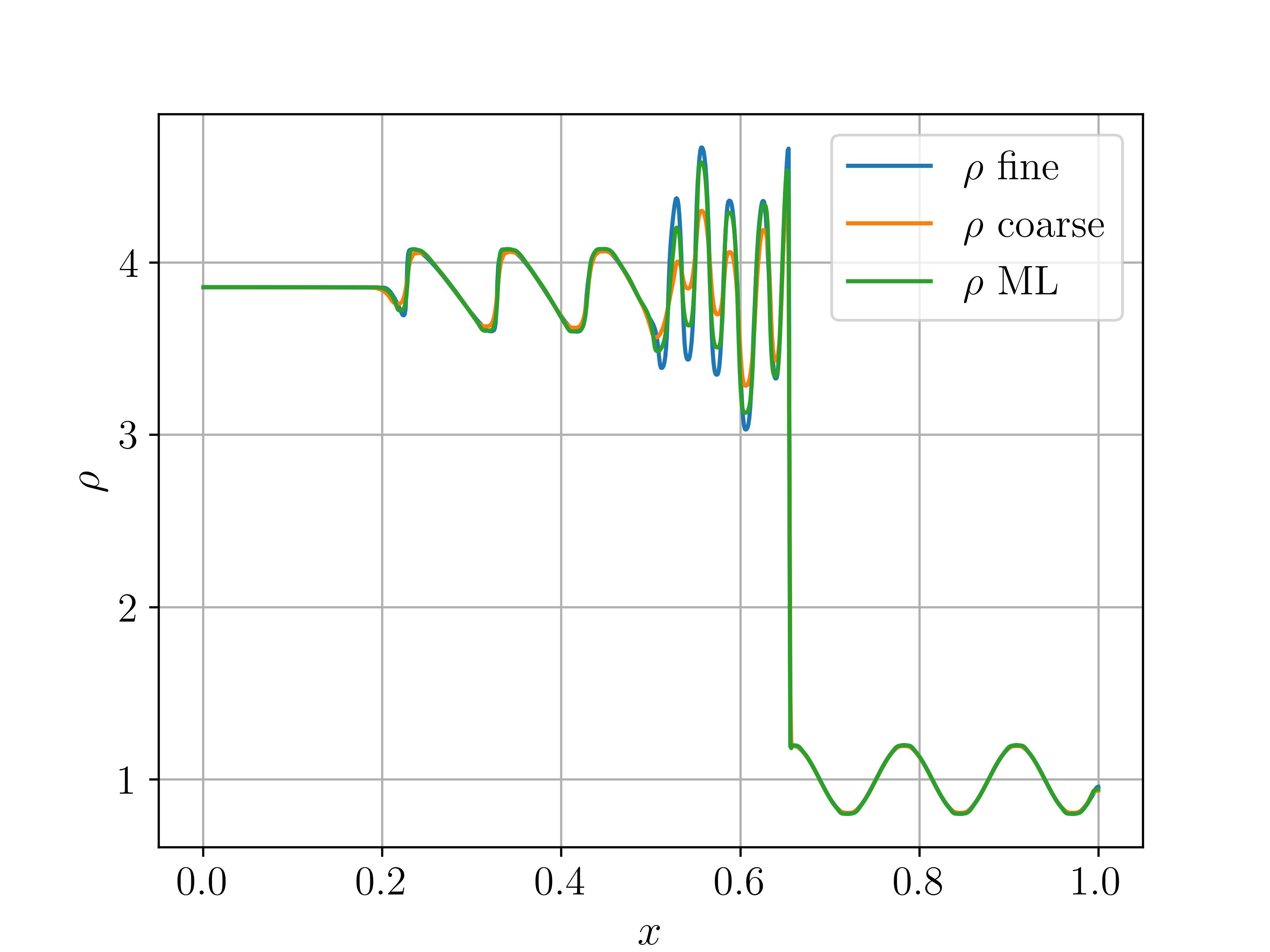}
         \caption{}
         \label{fig:sod_rho}
     \end{subfigure}
         
     \caption{Shu-Osher test case. 4096 cells for fine integration, 1024 cells for ML and for
coarse solutions, T = 0.156.}
\label{fig:shu_osher}
\end{figure}
Figures \ref{fig:sod} and \ref{fig:shu_osher} depict two reference test cases, showcasing how the ML solution compares to both the finite volume reference solution on the same grid and the reference solution on a finer grid. As seen on the Shu-Osher test case, the ML scheme has a great proficiency to well approximate sharp gradients. Besides, no difference can be seen on more regular part of the solution between the ML and the reference solutions. 
\par
On Figure \ref{fig:sod}, the error between the ML solution and the fine solution over time is plotted.The reference error between the coarse and the fine simulations is compared to the ML-based scheme error. We can conclude that on the Sod test case, excellent accuracy has been achieved. 

\subsection{2D Euler equations}
\par
\begin{figure}[H]
     \centering
     \begin{subfigure}[b]{0.3\textwidth}
         \centering
         \includegraphics[width=\textwidth]{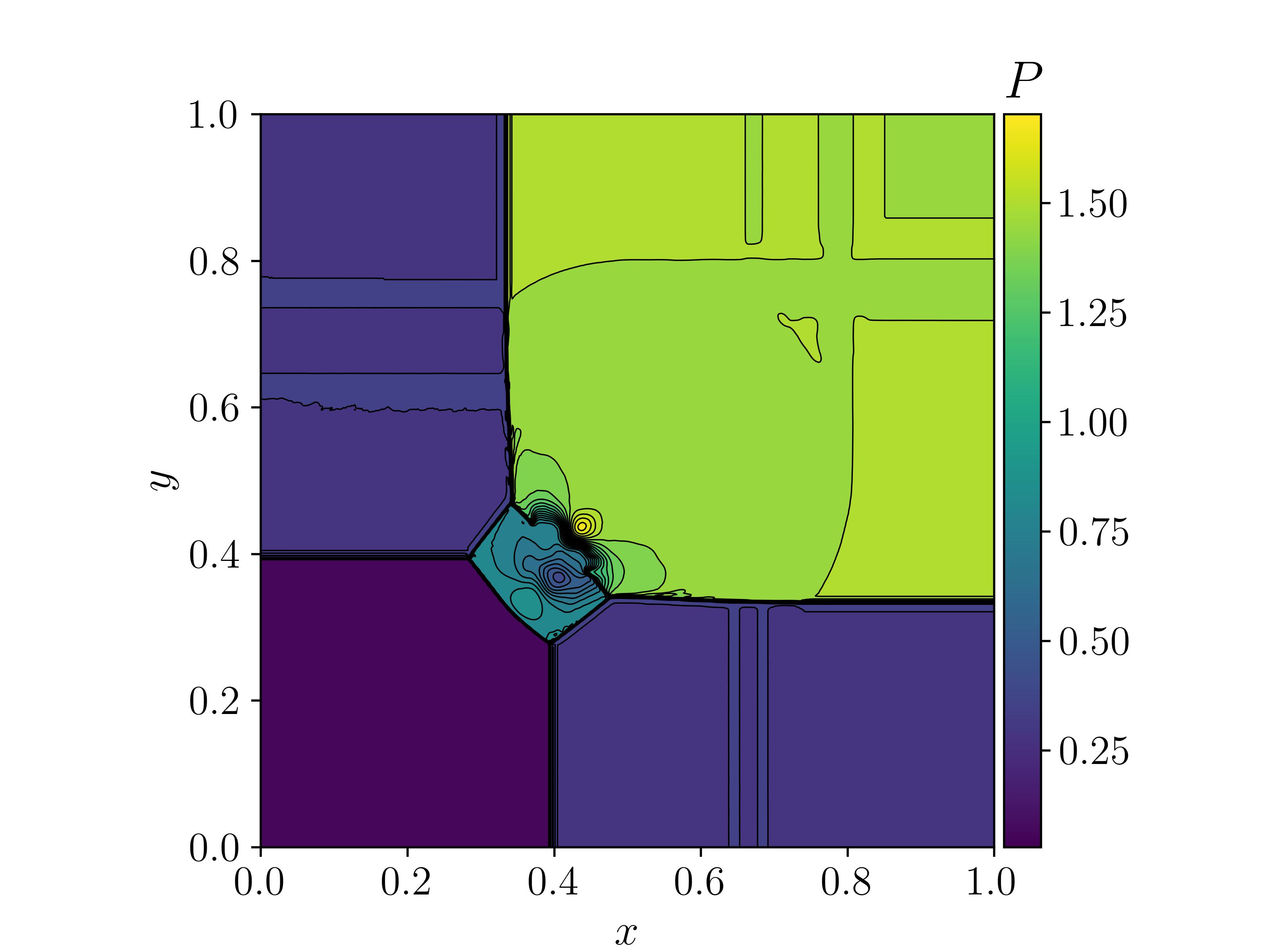}
         \caption{}
         \label{fig:sod_rho}
     \end{subfigure}
     \hfill
     \begin{subfigure}[b]{0.3\textwidth}
         \centering
         \includegraphics[width=\textwidth]{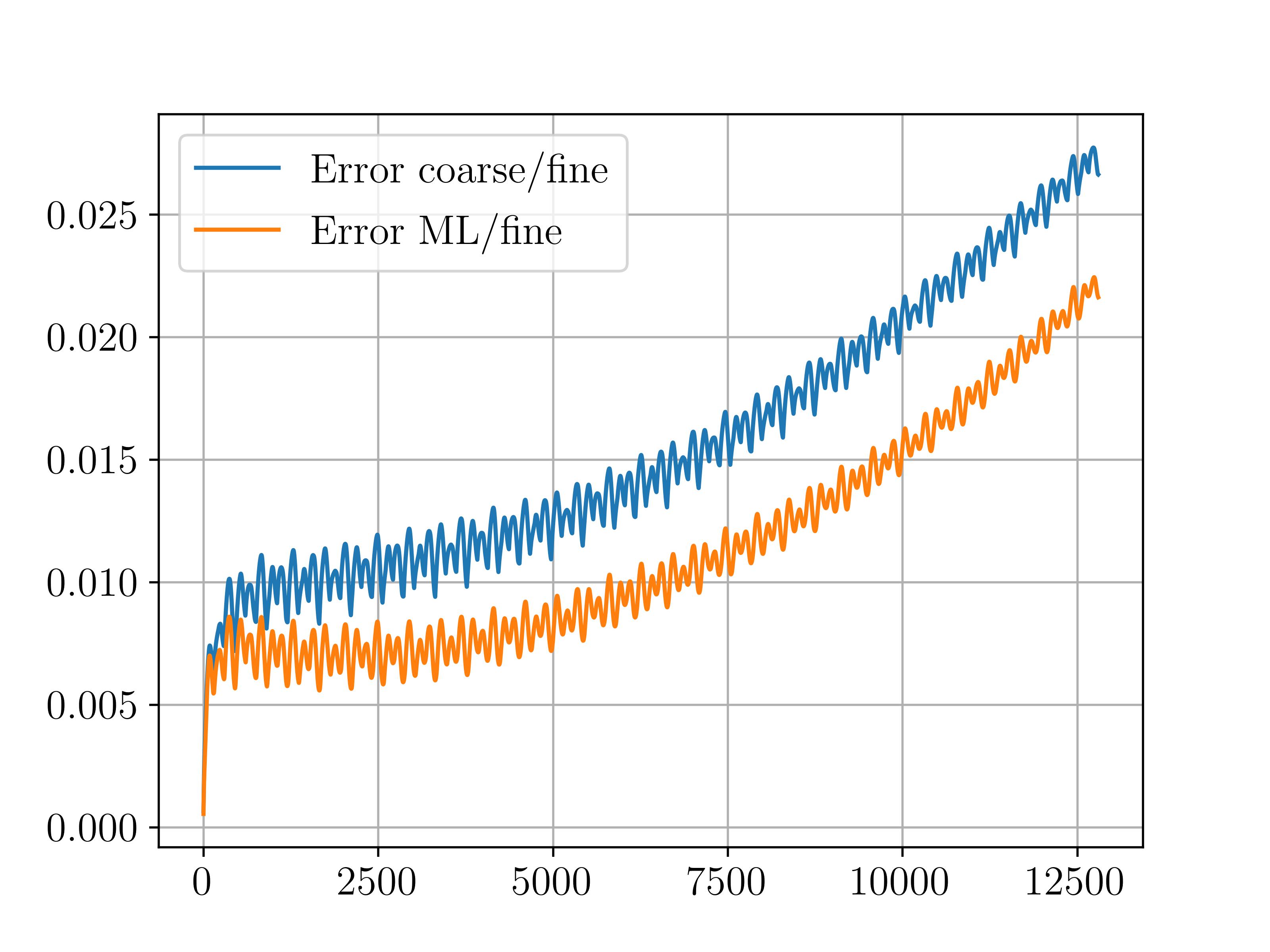}    
         \caption{}
         \label{fig:sod_rho_zoom}
     \end{subfigure}
     \caption{Results for the 2D Riemann problem case 3 \citep{liska2003comparison}. Density color map is overlayed by 30 density contours. The computations were performed on the  $(x, y) \in (0, 1) \times (0, 1)$ square on a 512x512 grid.}
     \label{fig:test_case_3}
\end{figure}

\begin{figure}[H]
     \centering
     \begin{subfigure}[b]{0.3\textwidth}
         \centering
         \includegraphics[width=\textwidth]{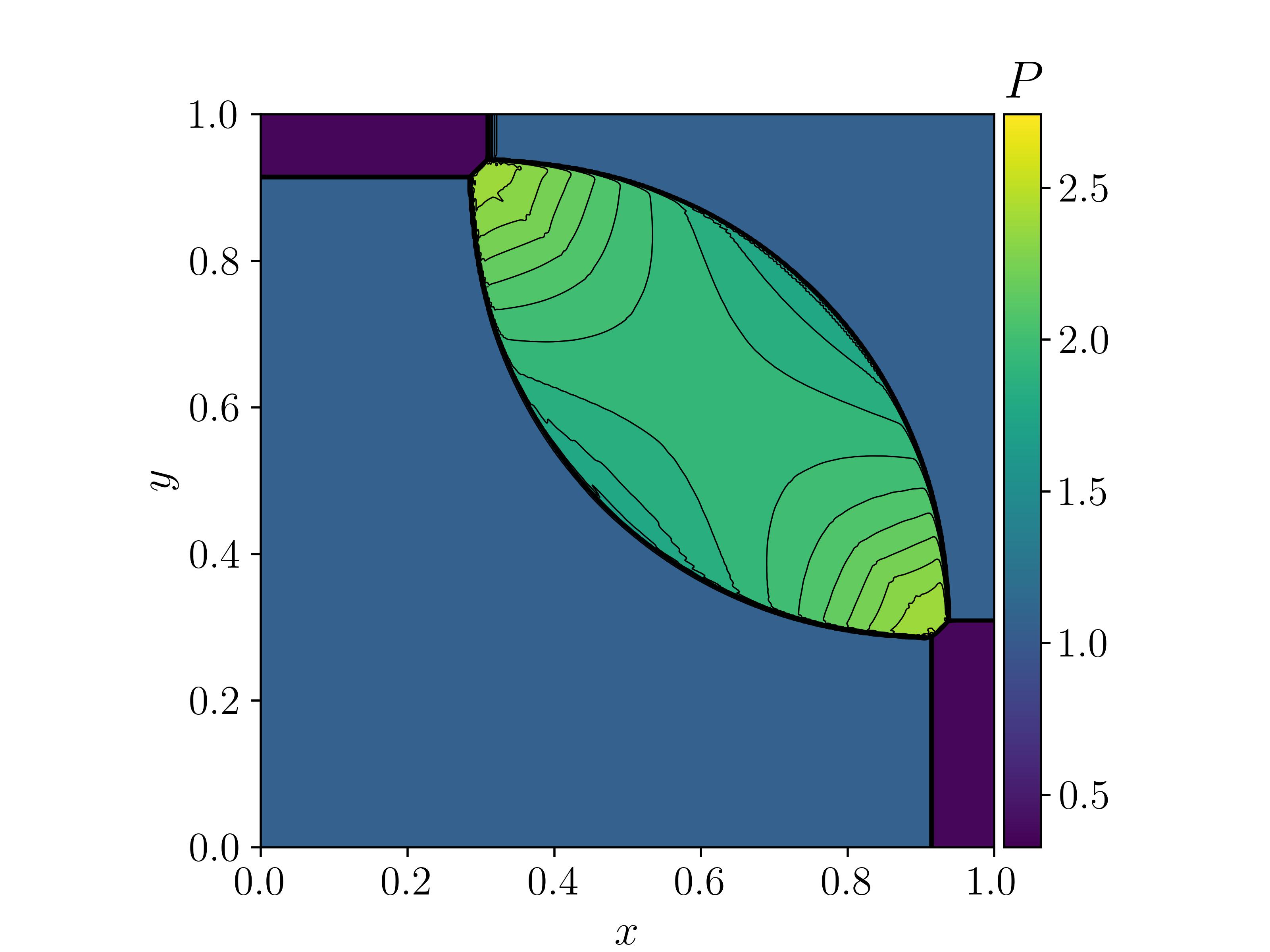}
         \caption{}
         \label{fig:sod_rho}
     \end{subfigure}
     \hfill
     \begin{subfigure}[b]{0.3\textwidth}
         \centering
         \includegraphics[width=\textwidth]{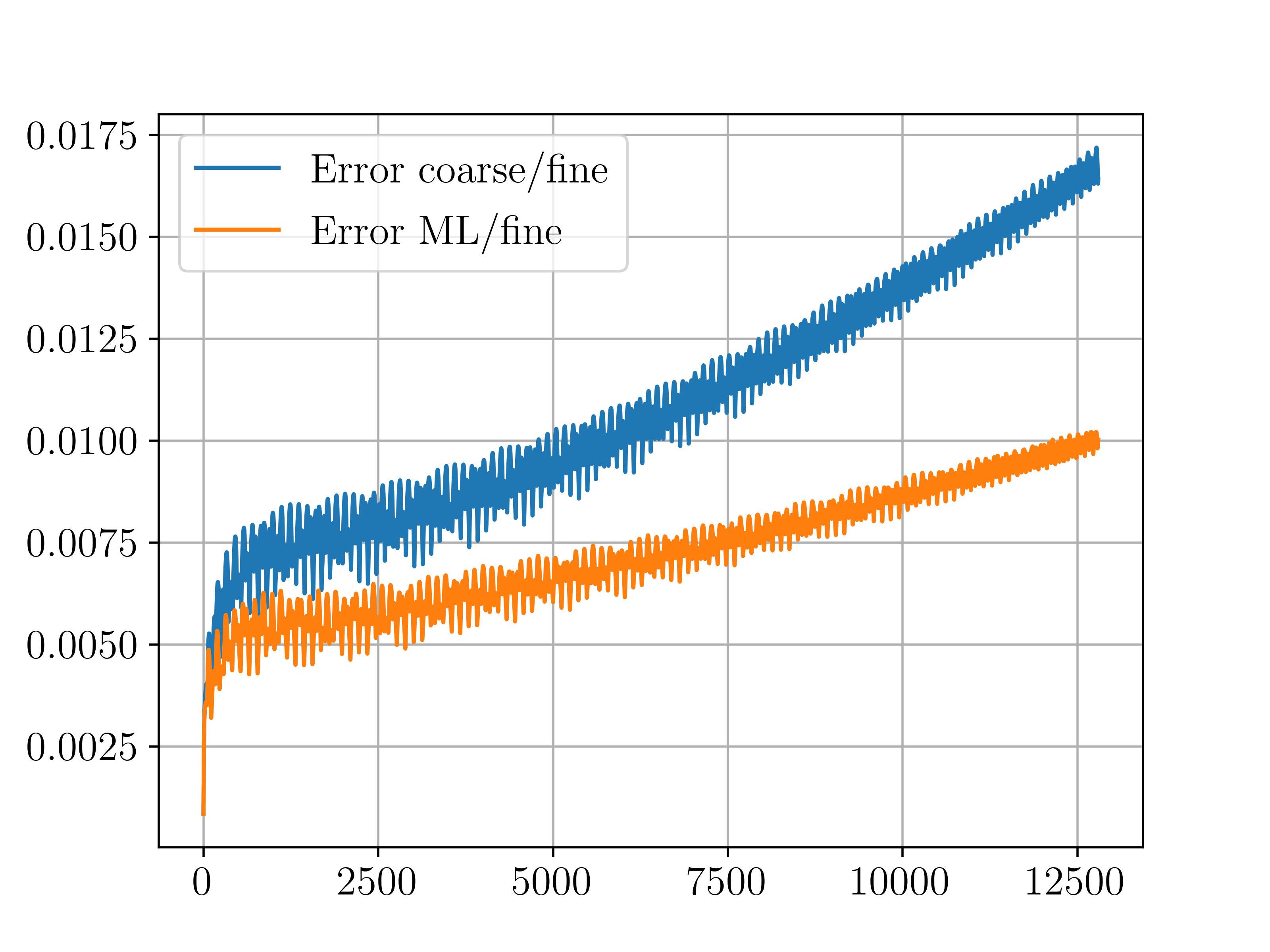}
         \caption{}
         \label{fig:sod_rho_zoom}
     \end{subfigure}
     \caption{Results for the 2D Riemann problem case 4 \citep{liska2003comparison}. Density color map is overlayed by 30 density contours. The computations were performed on the  $(x, y) \in (0, 1) \times (0, 1)$ square on a 512x512 grid.}
     \label{fig:test_case_4}
\end{figure}

\begin{figure}[H]
     \centering
     \begin{subfigure}[b]{0.3\textwidth}
         \centering
         \includegraphics[width=\textwidth]{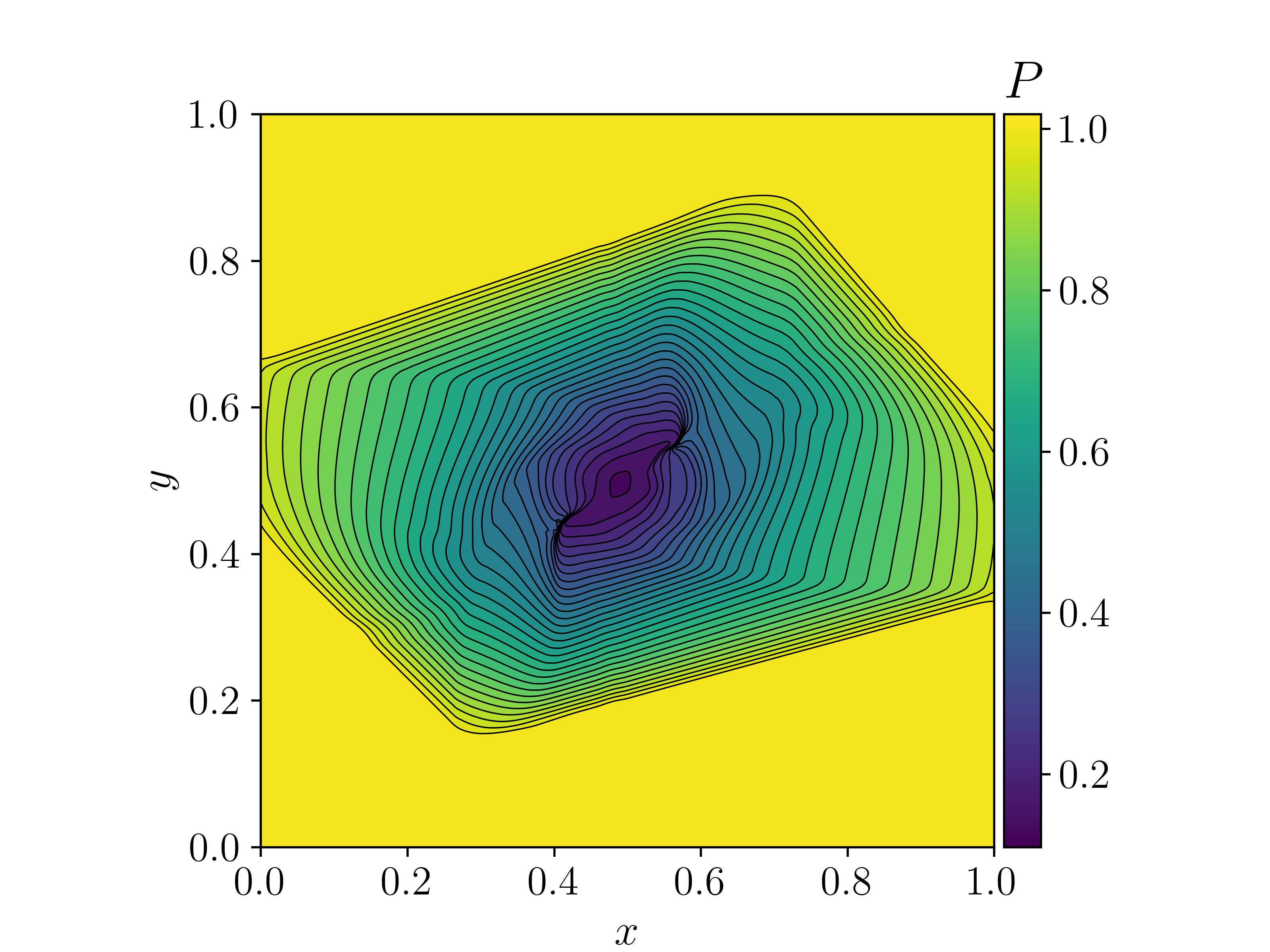}
         \caption{}
         \label{fig:sod_rho}
     \end{subfigure}
     \hfill
     \begin{subfigure}[b]{0.3\textwidth}
         \centering
         \includegraphics[width=\textwidth]{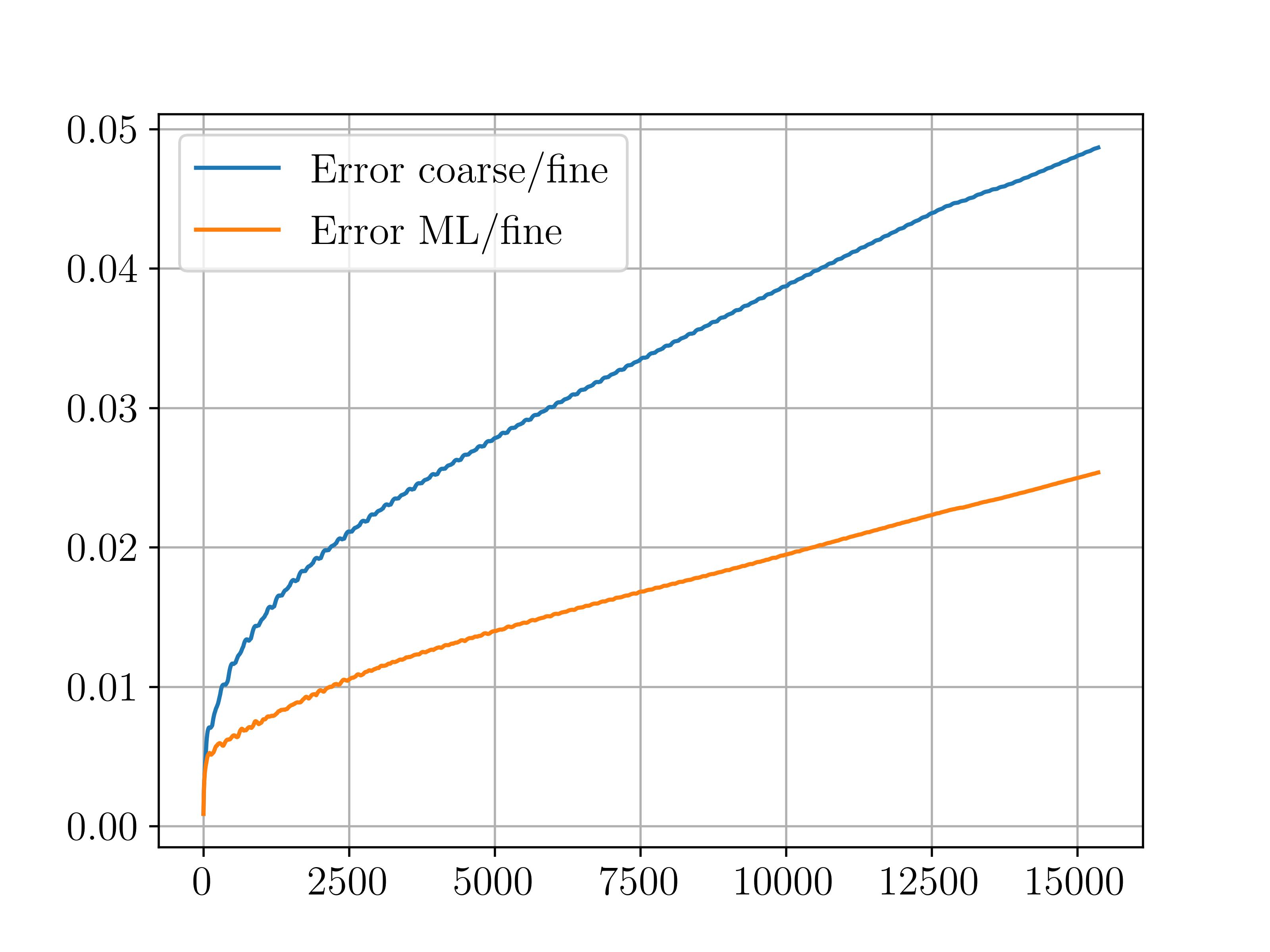}
         \caption{}
         \label{fig:sod_rho_zoom}
     \end{subfigure}
     \caption{Results for the 2D Riemann problem case 6 \citep{liska2003comparison}. Density color map is overlayed by 30 density contours. The computations were performed on the  $(x, y) \in (0, 1) \times (0, 1)$ square on a 512x512 grid.}
     \label{fig:test_case_6}
\end{figure}

\begin{figure}
     \centering
     \begin{subfigure}[b]{0.3\textwidth}
         \centering
         \includegraphics[width=\textwidth]{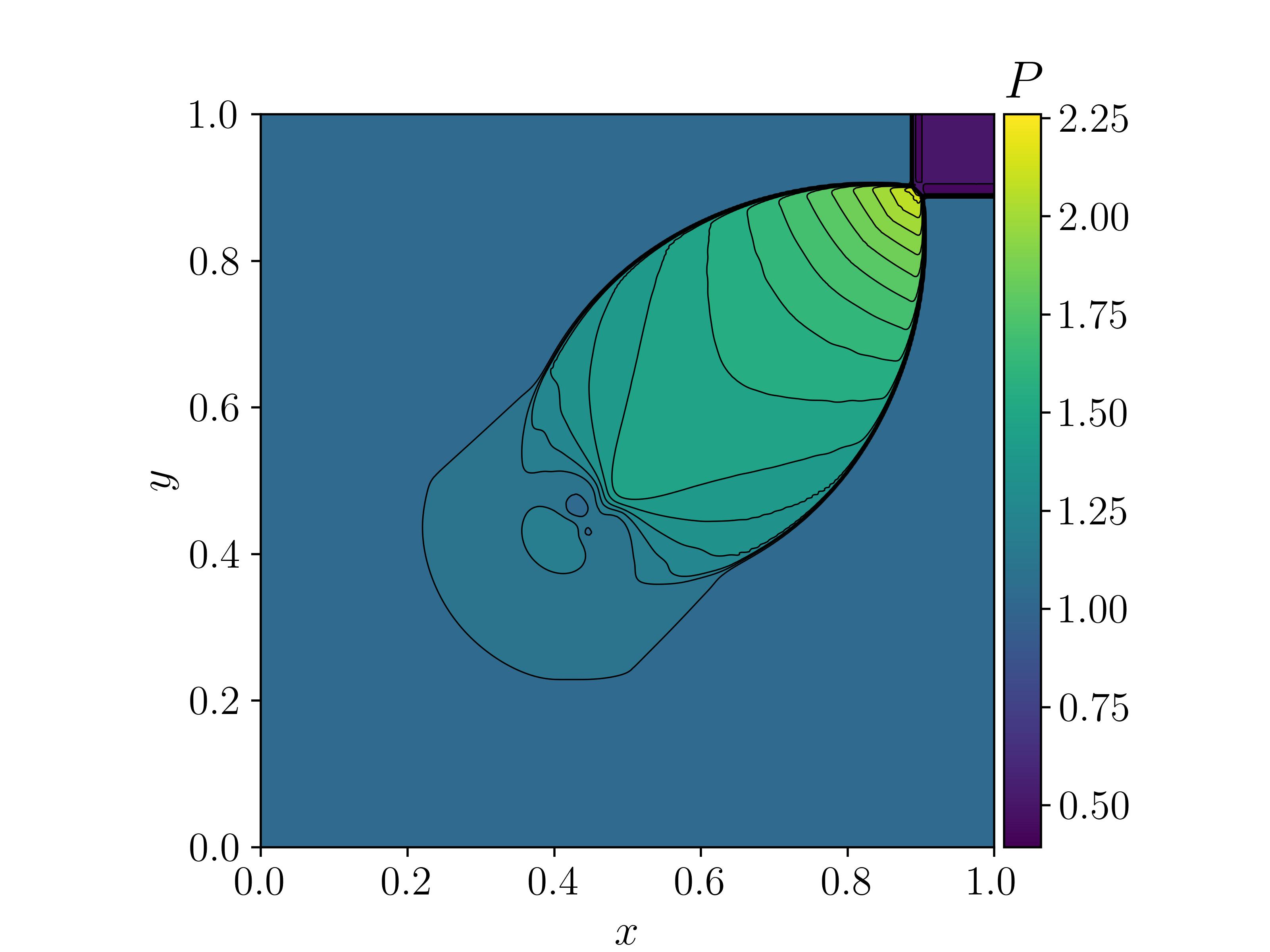}
         \caption{}
         \label{fig:sod_rho}
     \end{subfigure}
     \hfill
     \begin{subfigure}[b]{0.3\textwidth}
         \centering
         \includegraphics[width=\textwidth]{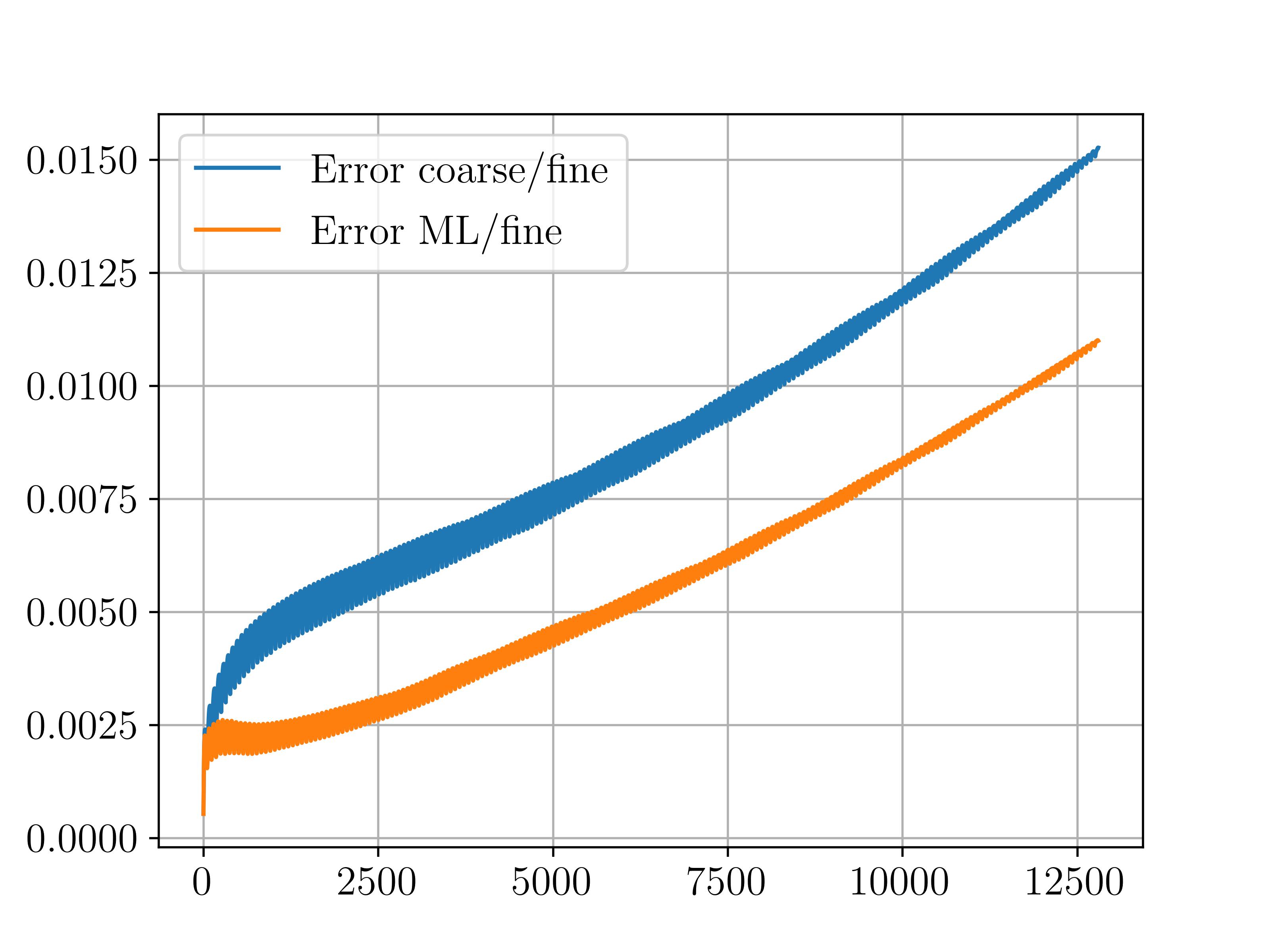}    
         \caption{}
         \label{fig:sod_rho_zoom}
     \end{subfigure}
     \caption{Results for the 2D Riemann problem case 12 \citep{liska2003comparison}. Density color map is overlayed by 30 density contours. The computations were performed on the  $(x, y) \in (0, 1) \times (0, 1)$ square on a 512x512 grid.}
     \label{fig:test_case_12}
\end{figure}

 
All test cases presented on Figures \ref{fig:test_case_3}, \ref{fig:test_case_4}, \ref{fig:test_case_6} and \ref{fig:test_case_12} are reference test cases for the 2D Euler equations presented by \cite{liska2003comparison}. Each figures shows the pressure in 2D and the error compared to a reference solution computed on a finer grid. Each test case is of different nature in order to test all kind shocks, rarefaction or contact-slip. All of the Riemann problems are proposed in such a way that the solutions of all four 1D Riemann problems between quadrants have exactly one wave (shock, rarefaction or contact-slip). 
\par 
The different test cases trials the robustness of the algorithm. Following \cite{lax1998solution} and defining $R$ for rarefaction, $S$ for shock and $J$ for contact-slip. Starting at the left side and going clockwise, the cases are: \textbf{Case 3}: $S, S, S, S$, \textbf{Case 4}: $S, S, S, S$, \textbf{Case 6}: $J, J, J, J$, \textbf{Case 12}: $J, S, S, J$.
\paragraph{Case 3} One can notice on Figure \ref{fig:test_case_3} the artefacts remaining on two segments of the initial discontinuities between the upper right quadrant and the upper left quadrant as well as on the lower right quadrants. Looking more closely at the intersection between the four quadrants, one can notice a lack of symmetry between the two plumes. 
\paragraph{Case 4} Compared to test case 3, this case preserves well the symmetry.

\paragraph{Case 6 } In this scenario, there are four instances of contact slip present across all four quadrants. This particular case stands out as one of the most effectively resolved. Additionally, the preservation of symmetry is notably high.

\paragraph{Case 12 } Presented Figure \ref{fig:test_case_12}, the key issue here is the resolution of the stationary contacts bordering the lower left quadrant. The data for this problem is symmetric about the $(0,0)-(1,1)$ diagonal and once again the scheme preserves well the symmetry of the case. 

\begin{figure}[H]
     \centering
     \begin{subfigure}[b]{0.3\textwidth}
         \centering
         \includegraphics[width=\textwidth]{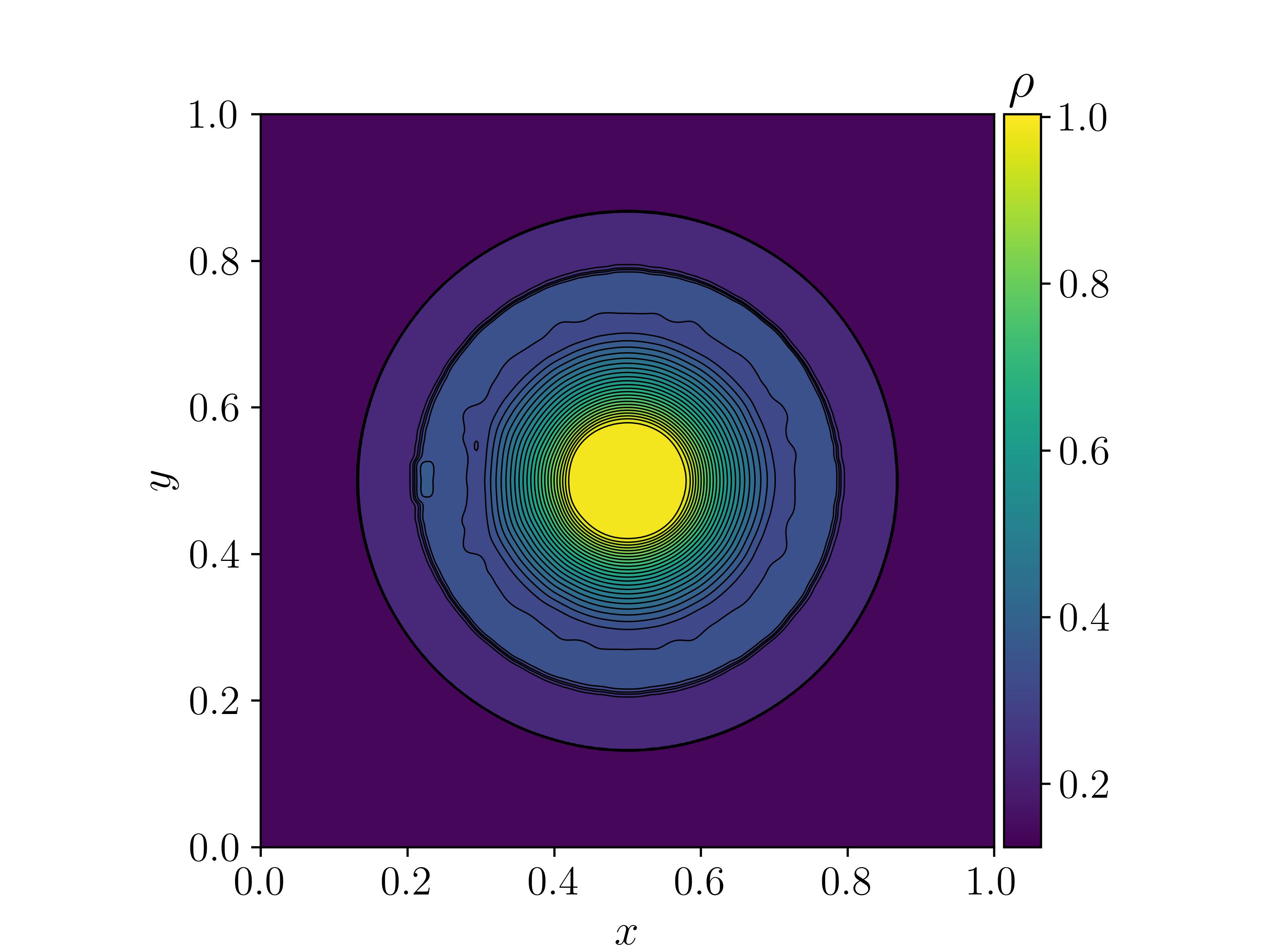}
         \caption{}
         \label{fig:sod_rho}
     \end{subfigure}
     \hfill
     \begin{subfigure}[b]{0.3\textwidth}
         \centering
         \includegraphics[width=\textwidth]{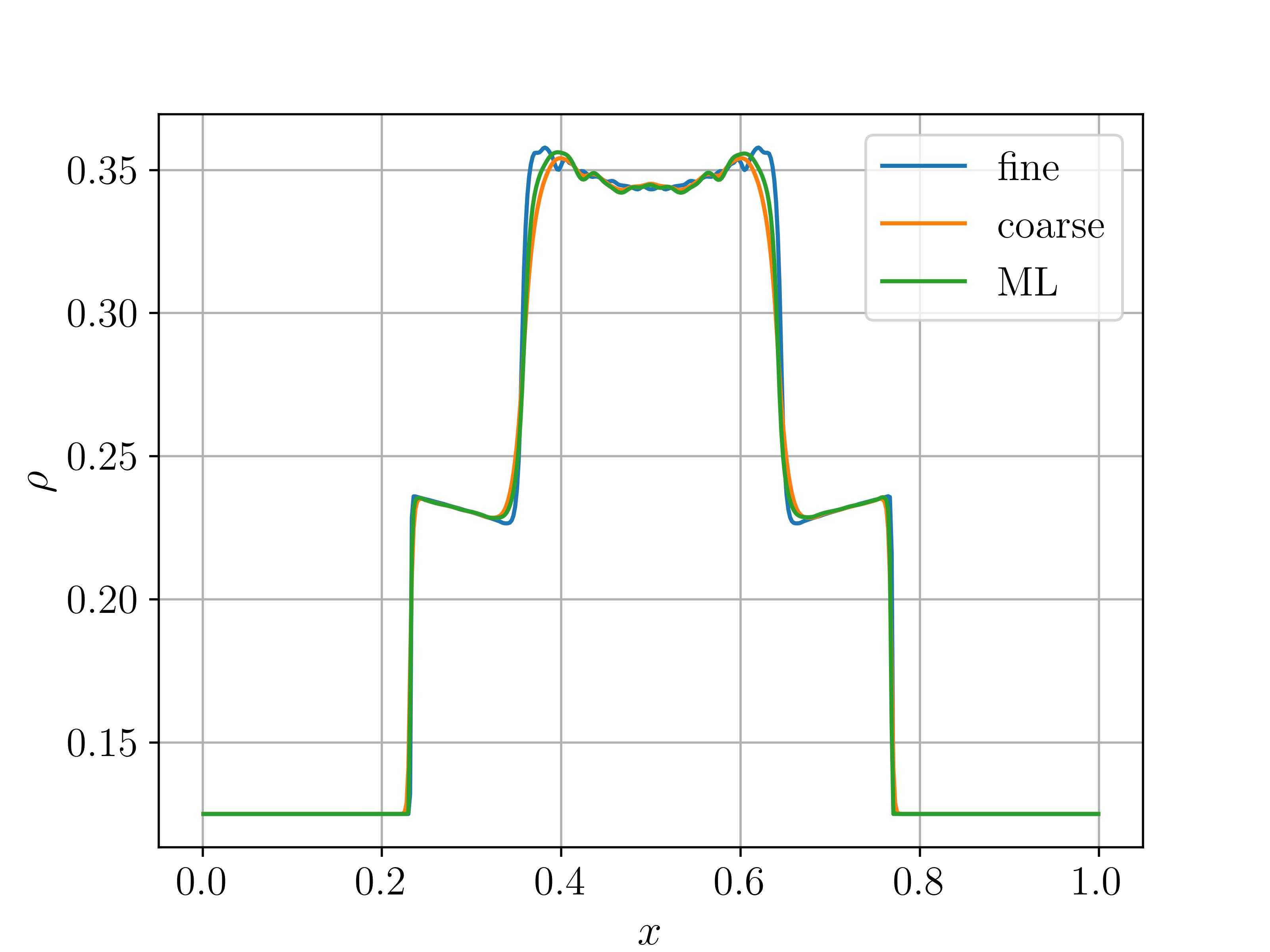}    
         \caption{}
         \label{fig:sod_rho_zoom}
     \end{subfigure}
     \caption{Results for the 2D explosion test case \citep{liska2003comparison}. 
     Density color map is overlayed by 30 density contours. 
     The computations were performed on the  $(x, y) \in (0, 1) \times (0, 1)$ square on a 256x256 grid. 
     The figure on the right represents a slice of the density along the $x$-axis at $y = 0.5$.
     }
     \label{fig:explosion}
\end{figure}

\paragraph{Explosion test case} The explosion problem proposed in \cite[sec. 17.1.1]{toro2013riemann} is a circularly symmetric 2D problem with initial circular region of higher density and higher pressure. In particular we set the center of the circle at $(0.5,0.5)$, its radius to $0.4$ and compute on a $(x, y) \in (0, 1.) \times (0, 1.)$ square. Density and pressure are $\rho_i = 1$, $p_i = 1$ inside the circle and $\rho_o = 0.125$, $p_o = 0.1$ outside. The test case presents much stronger shock and an instability which is well resolved here. The central part displayed on the right part of Figure \ref{fig:explosion} is well approximated but the outside instabilities are closer to the coarse reference solver. Noting that the reference solver also presents instabilities which the ML scheme seems to have dampened.

\paragraph{Forward facing step} The forward facing step test case described by \cite{woodward1984numerical} is a challenging test case mainly due to the shock interaction taking place at the corner of the step. This test case also has multiple shocks fronts across the domain and different boundary conditions (supersonic inflow/outflow, wall slip). This makes it an ideal test case for benchmarking. So, we treat it as a two-dimensional problem. Initially the wind tunnel is filled with a gamma-law gas, with $\gamma = 1.4$, with $(\rho, u , v, p ) = (1.4, 3., 0., 1.)$ as initial condition and a left-hand inlet boundary condition.
\begin{figure}[H]
     \centering
     \begin{subfigure}[b]{0.52\textwidth}
         \centering
         \includegraphics[width=\textwidth]{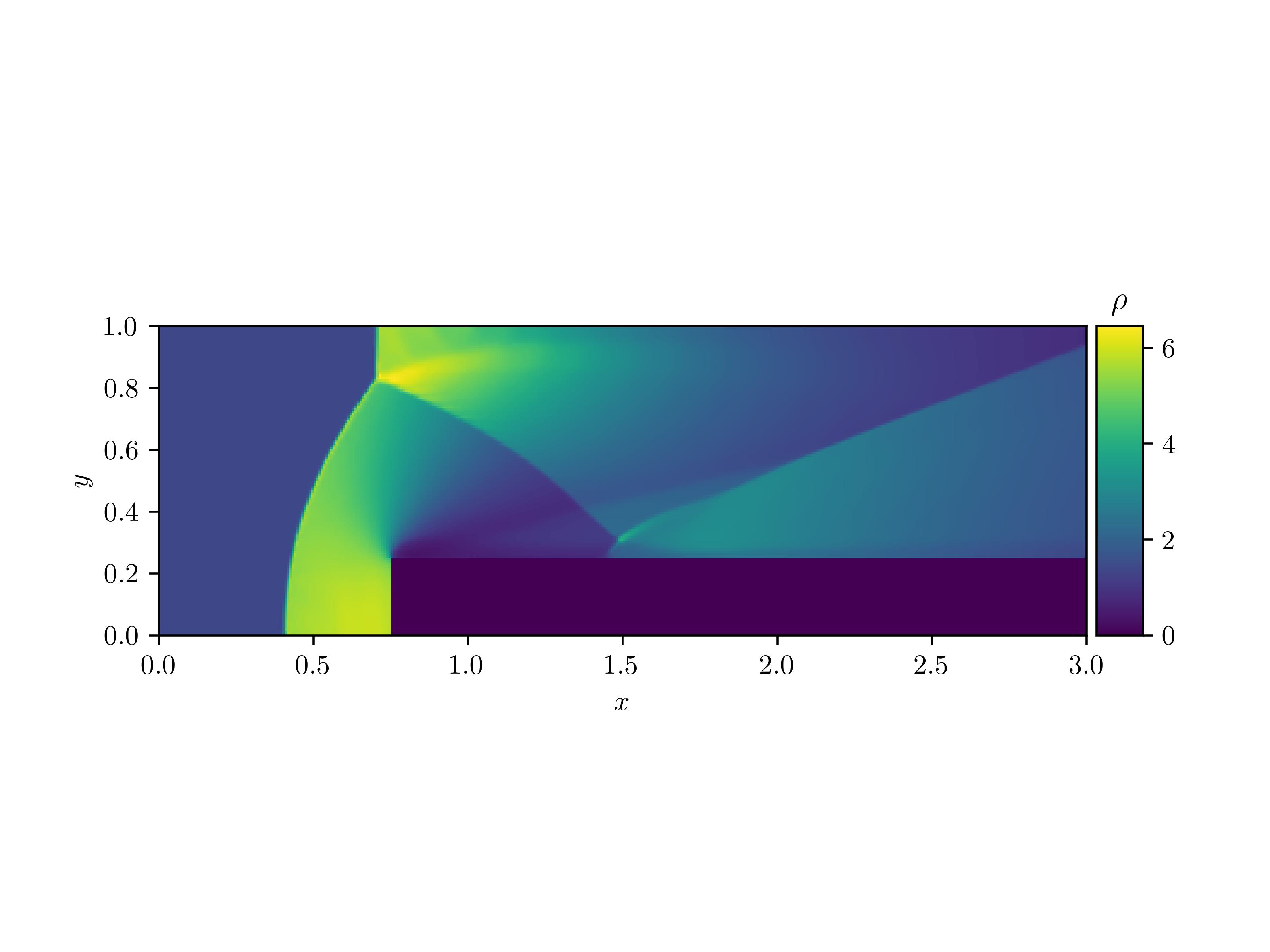}
         \caption{}
         \label{fig:sod_rho}
     \end{subfigure}
     \caption{Results for the 2D forward facing step \citep{woodward1984numerical}. Density color map is overlayed by 30 density contours. The computations were performed on a  $(x, y) \in (0, 3) \times (0, 1)$ rectangle on a 384x128 grid. 
     }
     \label{fig:forward_facing_step}
\end{figure}
\par 
The case is well resolved especially for all shocks fronts present in the solution. However as seen on the right of Figure \ref{fig:forward_facing_step}, an overshoot and some spurious oscillations are downgrading the quality of the the solution. This is in fact due to the double boundary condition at the corner. In spite of these oscillations, an $19\%$ increase in precision is observed on this test case for the same mesh. 

\subsection{Convergence study}
\begin{figure}[H]
     \centering
     \begin{subfigure}[b]{0.45\textwidth}
         \centering
         \includegraphics[width=\textwidth]{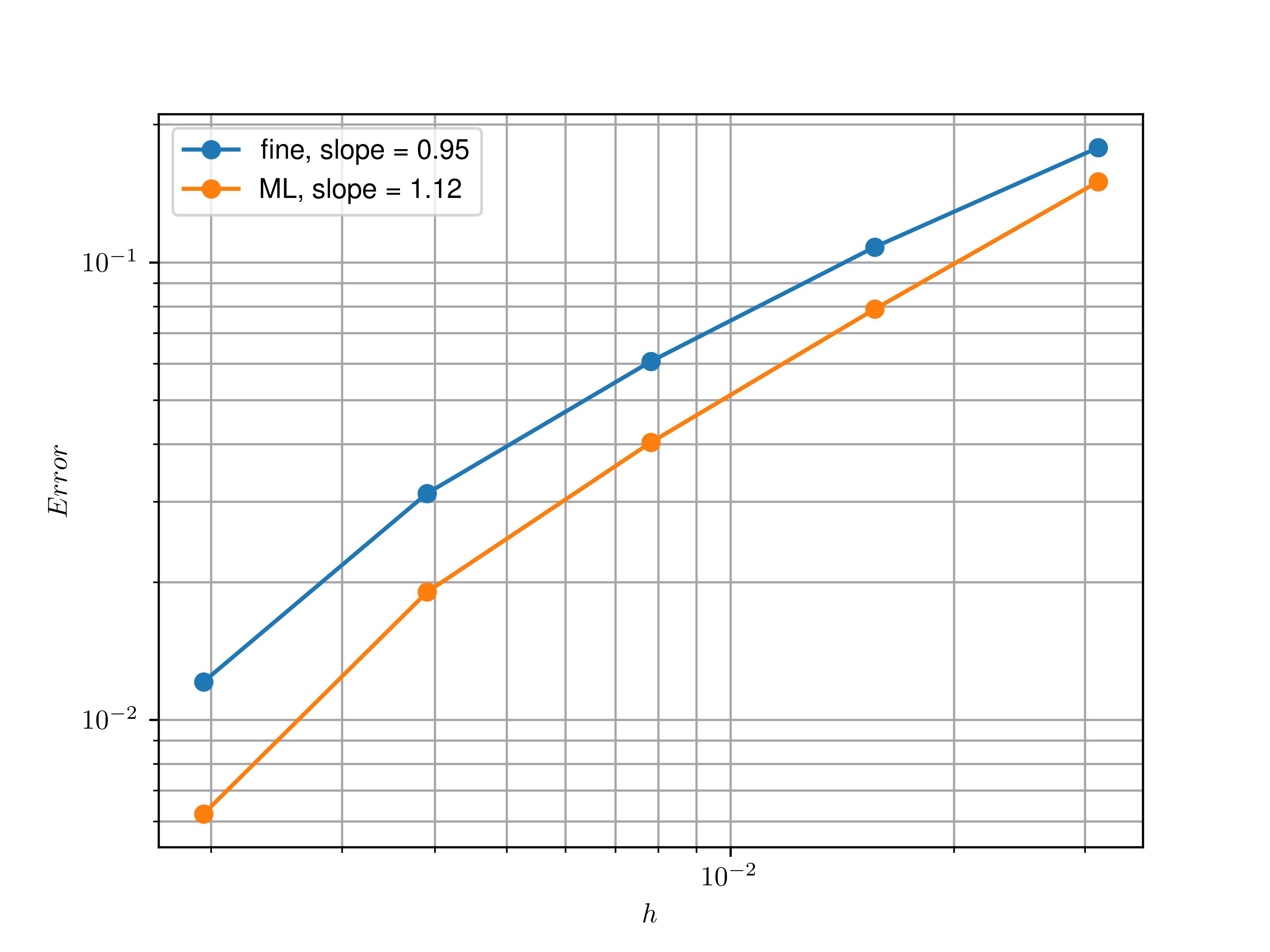}
         \caption{}
         \label{fig:convergence_Riemann}
     \end{subfigure}
     \hfill
     \begin{subfigure}[b]{0.45\textwidth}
         \centering
         \includegraphics[width=\textwidth]{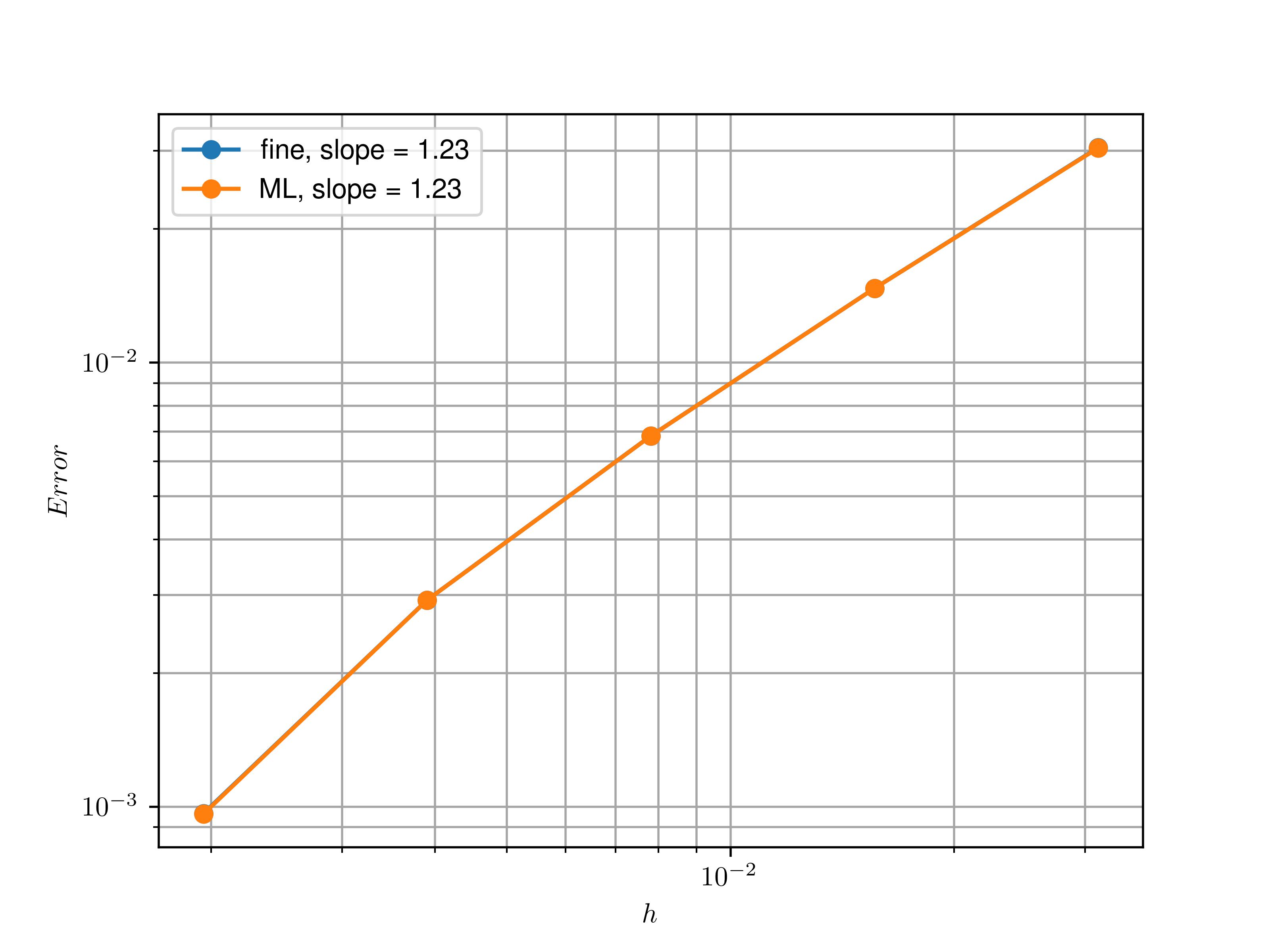}    
         \caption{}
         \label{fig:convergence_sinus}
     \end{subfigure}
     \caption{(a) Mesh convergence with Riemann initial conditions on a log-log scale (b) Mesh convergence with a sine-wave as initial conditions on a log-log scale }
     \label{fig:convergence}
\end{figure} 
Figure \ref{fig:convergence} presents the mesh convergence for the reference solver and the ML solver in 2D with two different initial conditions with periodic boundary conditions. The initial condition of Figure \ref{fig:convergence_Riemann} is a randomized Riemann-based initial condition as described for the dataset. The initial condition of Figure \ref{fig:convergence_sinus} is an advected sine-wave defined as $(\rho, u, v, p) = (sin(2\pi x), 1, 0, 1)$.
\par 
Figure \ref{fig:convergence} shows the utility of the machine learning method especially around shocks, indeed the methods with and without ML are identical for regular solution like for the sine-wave initial condition one. But in the case of multiple shocks, the ML method comes into effect compared to the reference finite volume solver meaning that the Machine Learning model method approximates well the derivatives around shocks.

\section{Computational performance}  \label{sec:Computational}
\begin{figure}[h]
    \includegraphics[width=0.48\textwidth]{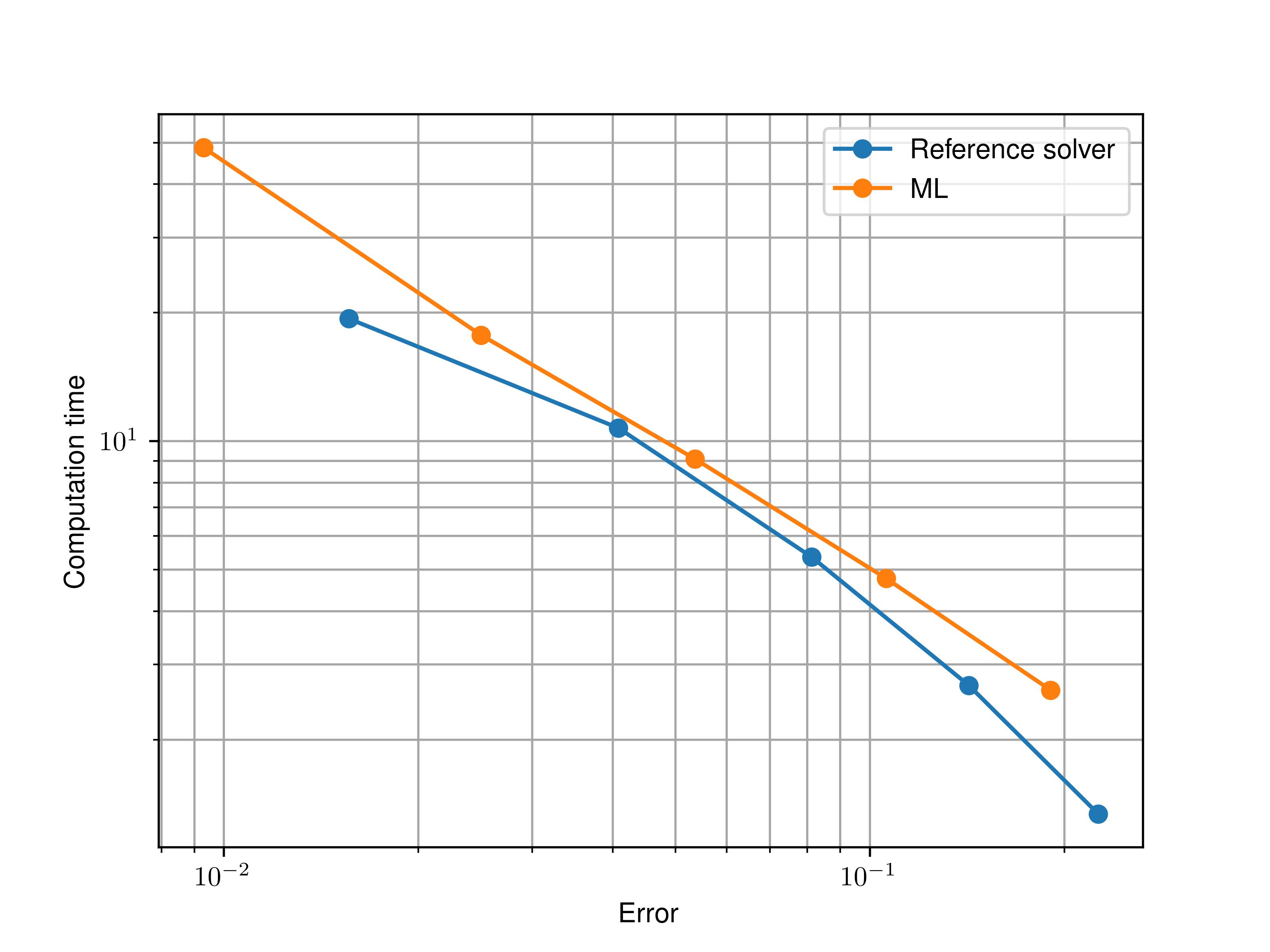}
    \caption{Computation time vs error on a log-log scale}
    \label{fig:computation}
\end{figure}
To test the computational performance of our Machine Learning method, during the convergence study the respective computation time is plotted along their respective error as seen on Figure \ref{fig:computation}. This method gives a good tendency of the algorithms by comparing the two main objectives in CFD : the computation time and the error.
\par 
One downside of the built-in  algorithm is the computation time. 
All calculations were done on a GPU Nvidia V100. As seen on Figure \ref{fig:computation}, for a same error the reference solver takes less computation time, especially for small errors or fine grids. 
\section{Discussion and conclusion}

We have adapted the method presented by \cite{bar2019learning} to hyperbolic partial differential equations (PDEs) and achieved notable enhancements in accuracy, robustness, and solution quality when solving conservation laws with neural networks. 
A number of regularizers have been introduced to ensure that the neural network solution to the forward problems is physically achievable, particularly in situations involving strong shocks. We also note that thanks to the limiter, the training is not done on unrolled steps, speeding-up notably the training time.
The algorithm has been studied on a linear equation to show theoretical results and gains. The method has also been tested on different equations and on various cases from the literature showing very good results in terms of accuracy and stability. 

It should be noted that the method is not without limitations. One identified drawback of the approach developed here is the increased computation time in comparison to the reference finite volume scheme. It would be beneficial to dedicate future studies to improving computation time, which is closely linked to the use of the TensorFlow library. Additionally, our tests have focused on Cartesian meshes. However, more flexible convolutions \citep{sanchez2020learning} could be employed to utilize the presented method on more flexible mesh structures with irregular discretizations. Furthermore, even regular CNNs can be adapted to account for non-uniform and stretched meshes \citep{chen2021towards}.

It bears noting that this paper represents a proof of concept. The potential for data-driven discretizations to facilitate the efficient resolution of previously intractable industrial problems is a promising avenue for future exploration.

\section{Acknowledgments}
We thank Piotr Mirowski from DeepMind for helpful discussions. 

\bibliography{bibliography}

\newpage
\appendix
\section*{Appendices}
\addcontentsline{toc}{section}{Appendices}
\renewcommand{\thesubsection}{\Alph{subsection}}
\subsection{Dataset}
\label{sec:dataset}
The training data set is generated using pre-defined parametrized functions as initial conditions on the domain $[0,1]$ in 1D or $[0,1]^2$ in 2D.
\par 
Periodic and non-slip boundary dataset have the same number of initial conditions picked at the same share at random and integrated on the same number of time steps.
\par
We define $Rect(x_0, x_1)$ the rectangular function in 1D defined as : 
\begin{equation*}
    Rect(x_0, x_1) = 
        \begin{cases}
      1 & \text{if $x\geq x_0$ and $x \leq x_1$}\\
      0 & \text{else}
    \end{cases}   
\end{equation*}
\subsubsection{Burgers}
In this section are presented the initial conditions used for training the model for the Burgers equation.
\begin{align*}
    f(x) &= (1 -\frac{1}{2} Rect(0.15, 0.35))\displaystyle \sum_{i=1}^N a_i \dfrac{sin(2\pi l_i x + \phi_i)}{3}\
\end{align*}
\par
With $\forall i, a_i\in [-0.5, 0.5], \phi_i\in [0, 2\pi]$ and $l_i \in \{4,...,N\}$. N is usually set to 20. 
\par
The training data is composed of 50 initial conditions integrated on 10,000 time steps making 490,000 samples. The test data is composed of only 5 different initial conditions.
\begin{figure}[h]
\captionsetup{justification=centering}
\centering
    \includegraphics[width=0.35\textwidth]{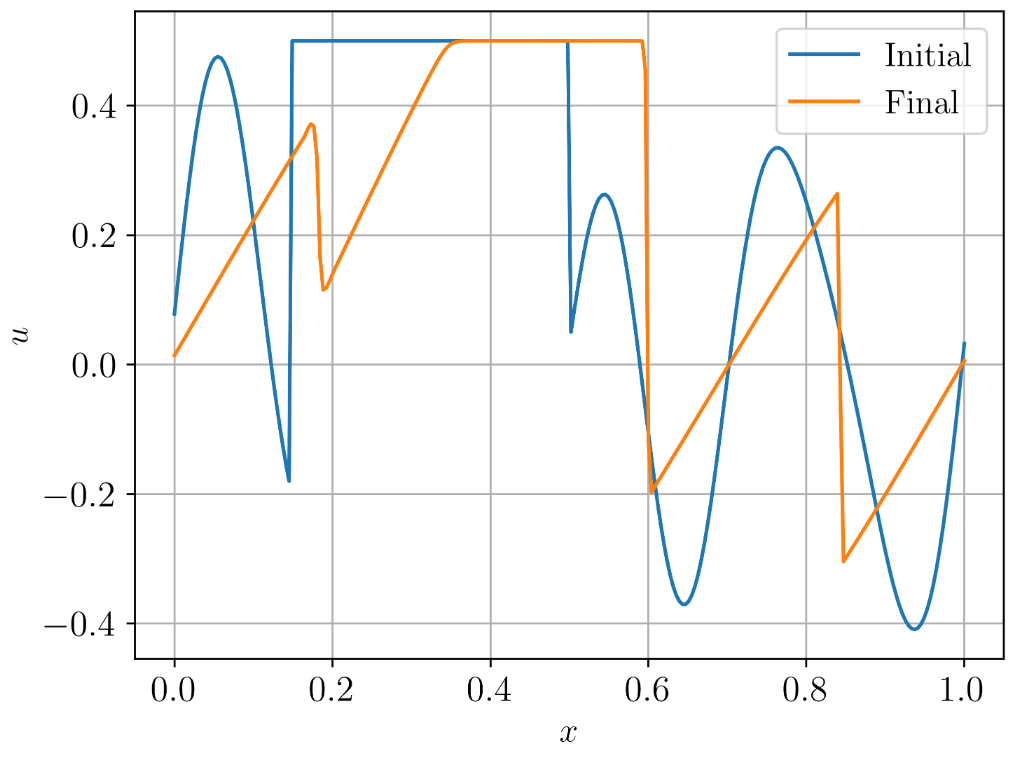}
    \caption{Initial and final solutions for the 1D Burger database with periodic boundary conditions.}
\label{fig:database_1D_burger}
\end{figure}
\subsubsection{Euler 1D}
In this section are presented the initial conditions used for training the model for the 1D Euler equation.
\begin{align*}
\textbf{f}_1(x) &=
    \begin{cases}
        f_{\rho} &= sin(2 \pi x + \phi_0\pi) + 1.2 + h_0\\ 
        f_{u} &= sin(2 \pi x + \phi_1\pi) + 1 + h_1\\
        f_{p} &= sin(2 \pi x + \phi_2\pi) + 1 + h_2
    \end{cases}\\
\textbf{f}_2(x) &=
    \begin{cases}
        f_{\rho} &= (a_0 - 0.5) Rect(x_0, 1) + h_0 + 0.7\\ 
        f_{u} &= a_1 Rect(x_0, 1) \\
        f_{p} &= -(a_2-0.5) Rect(x_0, 1) + 2 h_2 + 0.7  
    \end{cases}\\
\textbf{f}_3(x) &=
    \begin{cases}
        f_{\rho} &= a_0 Rect(x_0, x_1) + h_0 + 0.1\\ 
        f_{u} &= a_1 Rect(x_0, x_1) + h_1\\
        f_{p} &= a_2 Rect(x_0, x_1) + h_2
    \end{cases}
\end{align*}
All parameters are random parameters following a distribution $\mathcal{U}(0,1)$. We just set $0.2\leq|x_0 - x_1|\leq0.8$.
\par 
The training data is composed of 15 initial conditions composed of 60\% $\textbf{f}_1$, 20\% $\textbf{f}_2$ and 20\% $\textbf{f}_3$  integrated on 6,001 time steps making 90,000 samples. The test data is composed of only 5 different initial conditions also integrated on 6,001 time steps.
\begin{figure}[h]
\captionsetup{justification=centering}
\centering
    \begin{subfigure}[h]{0.35\textwidth}
        \includegraphics[width=\textwidth]{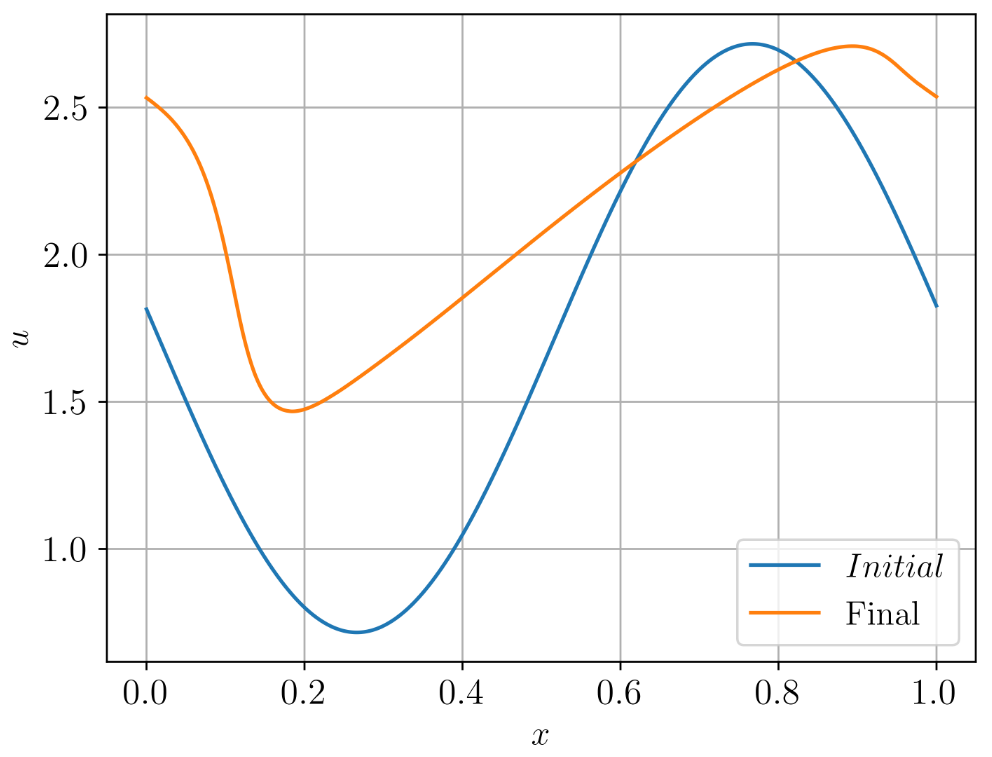}
        \caption{A representation of $\textbf{f}_1$}
    \end{subfigure}
    \begin{subfigure}[h]{0.35\textwidth}
        \includegraphics[width=\textwidth]{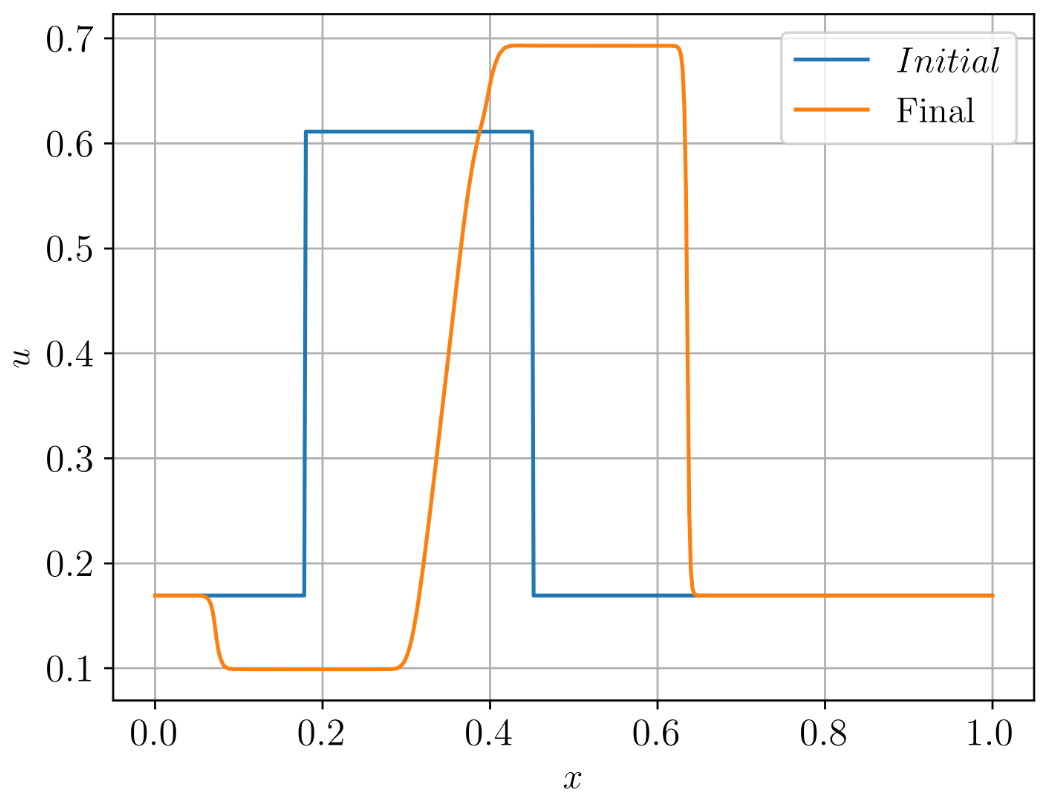}
        \caption{A representation of $\textbf{f}_3$}
    \end{subfigure}
\caption{Initial and final solutions for the 1D Euler database with periodic boundary conditions.}
\label{fig:database_1D_euler}
\end{figure}
\subsubsection{Euler 2D}
In this section are presented the initial conditions used for training the model for the 2D Euler equation.
\par
We first define the function 
\begin{equation*}
    \mathcal{C}(\textbf{x},\mathbf{p}) =
    \begin{cases}
        p_0 & \text{if $||(\textbf{x} - (0.5, 0.5)||_ 2\leq 0.125$}\\
        p_1 & \text{else if $x<0.5$ and $y<0.5$} \\
        p_2 & \text{else if $x\geq0.5$ and $y<0.5$} \\
        p_3 & \text{else if $x<0.5$ and $y\geq0.5$} \\
        p_4 & \text{else}
    \end{cases}
\end{equation*}
and the function
\begin{equation*}
    \mathcal{R}(\textbf{x},\mathbf{p}) =
    \begin{cases}
        p_0 & \text{if $||(\textbf{x} - (0.5, 0.5)||_ 1\leq 0.2$}\\
        p_1 & \text{else} 
    \end{cases}
\end{equation*}
with $\textbf{p}$ being parameters following a distribution $\mathcal{U}(0,1)$. 
\begin{align*}
\textbf{f}_1(\textbf{x}) &=
    \begin{cases}
        f_{\rho} &= a_0 sin(4 \pi x + \phi_0\pi) + a_1 sin(4 \pi y + \phi_1\pi) + 2 \\ 
        f_{u} &= 2(a_2 - 0.5) sin(4 \pi x + \phi_2\pi) + 2(a_3 - 0.5) sin(4 \pi y + \phi_3\pi)\\
        f_{v} &= 2(a_4 - 0.5) sin(4 \pi x + \phi_4\pi) + 2(a_5 - 0.5) sin(4 \pi y + \phi_5\pi)\\
        f_{p} &= a_6 sin(4 \pi x + \phi_6\pi) + a_6 sin(4 \pi y + \phi_7\pi) + 2
    \end{cases} \\
\textbf{f}_2(\textbf{x}) &=
    \begin{cases}
        f_{\rho} &= \mathcal{C}(\textbf{x},\frac{1}{2}\mathbf{p}_0 + \frac{1}{2}) \\ 
        f_{u} &= \mathcal{C}(\textbf{x},2(\mathbf{p}_1-1))\\
        f_{v} &= \mathcal{C}(\textbf{x},2(\mathbf{p}_2-1))\\
        f_{p} &= \mathcal{C}(\textbf{x},0.8\mathbf{p}_3 + 0.2)
    \end{cases} \\
\textbf{f}_3(\textbf{x}) &=
    \begin{cases}
        f_{\rho} &= \mathcal{R}(\textbf{x},\mathbf{p}_0 + \frac{1}{2}) \\ 
        f_{u} &= \mathcal{R}(\textbf{x},2(\mathbf{p}_1-1))\\
        f_{v} &= \mathcal{R}(\textbf{x},2(\mathbf{p}_2-1))\\
        f_{p} &= \mathcal{R}(\textbf{x},\mathbf{p}_3 + \frac{1}{2})
    \end{cases}
\end{align*}
\par
The training data is composed of 10 initial conditions composed of 37.5\% $\textbf{f}_1$, 37.5\% $\textbf{f}_2$ and 25\% $\textbf{f}_3$  integrated on 6,001 time steps making 60,000 samples. The test data is composed of only 2 different initial conditions also integrated on 6,001 time steps.

\begin{figure}[h]
\captionsetup{justification=centering}
    \begin{subfigure}[h]{0.45\textwidth}
        \includegraphics[width=\textwidth]{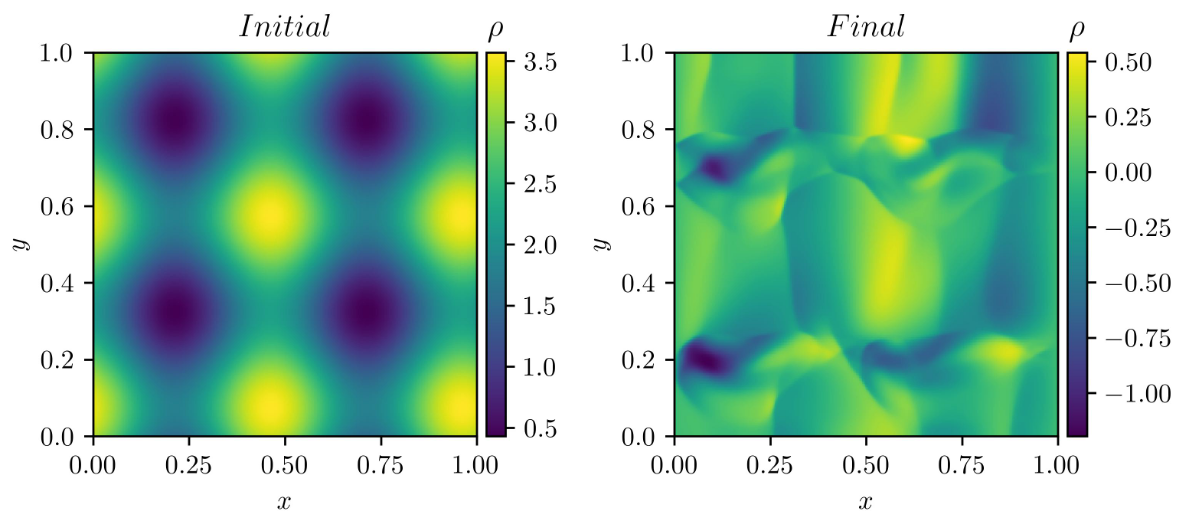}
        \caption{A representation of $\textbf{f}_1$}
    \end{subfigure}
    \begin{subfigure}[h]{0.45\textwidth}
        \includegraphics[width=\textwidth]{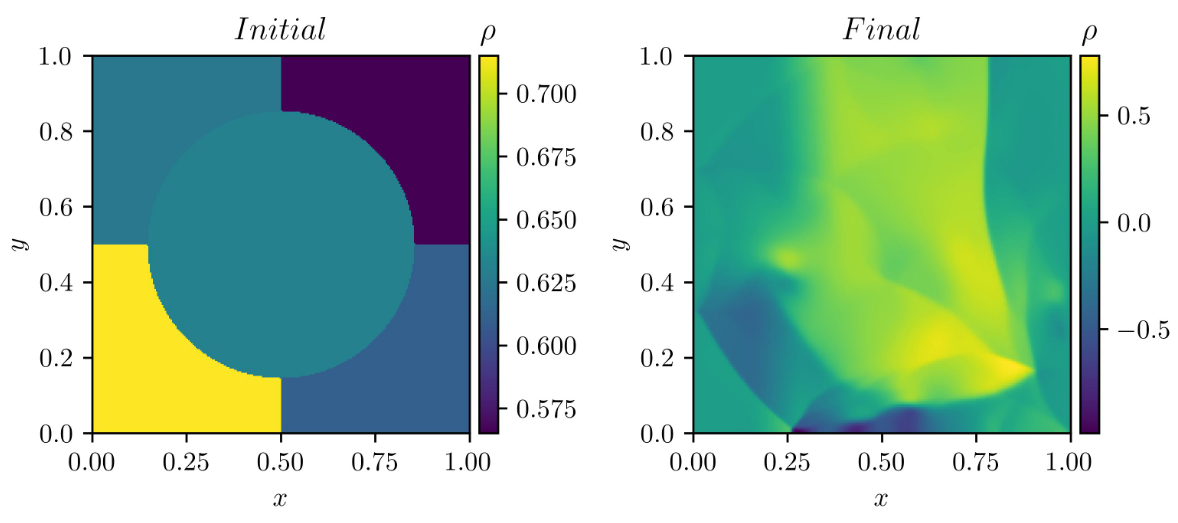}
        \caption{A representation of $\textbf{f}_2$}
    \end{subfigure}
    \begin{subfigure}[h]{0.45\textwidth}
        \includegraphics[width=\textwidth]{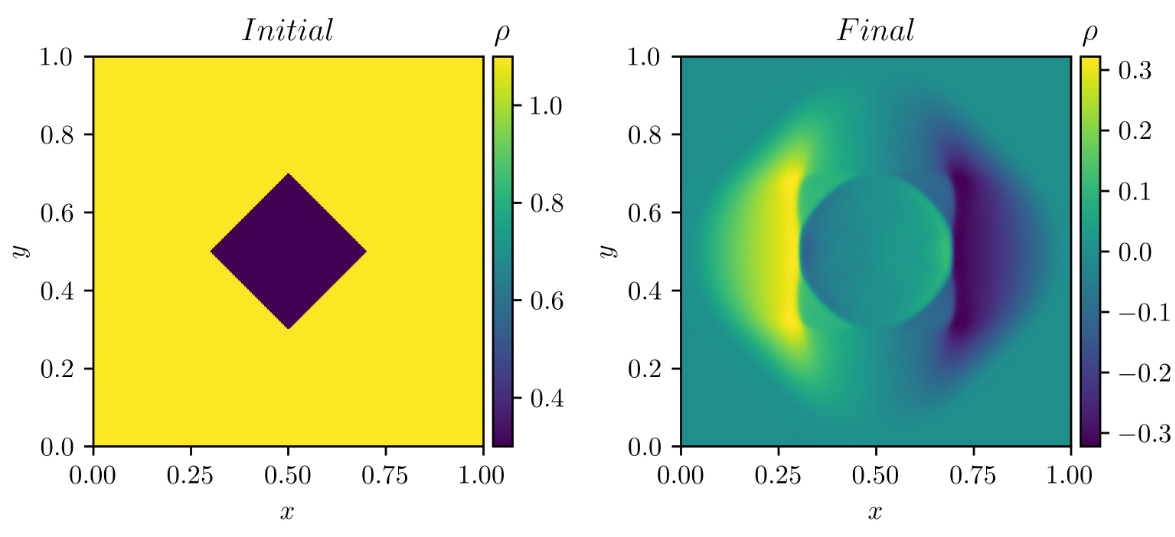}
        \caption{A representation of $\textbf{f}_3$}
    \end{subfigure}
\caption{Initial and final solutions for the 2D Euler database with periodic boundary conditions.}
\label{fig:database_2D_euler}
\end{figure}
\end{document}